\documentclass[a4, 11pt.]{amsart}
\usepackage{amssymb,amscd, hyperref, bm, color }
\usepackage{amsmath}
\usepackage{enumerate}
\usepackage{stmaryrd}  
\usepackage[mathcal]{eucal}
\usepackage{array,float}
\usepackage{longtable}
\usepackage{multirow}
\usepackage{cite}
\usepackage[left=3cm, right=3cm]{geometry}

\usepackage{xy}
\input xy
\xyoption{all}
\usepackage{pdflscape} 
\usepackage{hyperref}
\hypersetup{
	colorlinks=blue,
	linkcolor=blue,
	filecolor=magenta,      
	urlcolor=cyan,
}

\usepackage[draft]{graphicx}

\numberwithin{equation}{section}
\newtheorem{thm}{Theorem}[section]
\newtheorem{proposition}[thm]{Proposition}

\newtheorem{corollary}[thm]{Corollary}
\newtheorem{lemma}[thm]{Lemma}

\theoremstyle{definition}
\newtheorem{definition}[thm]{Definition}
\newtheorem*{remark*}{Remark}
\newtheorem{remark}[thm]{Remark}

\newcommand{\compcent}[1]{\vcenter{\hbox{$#1\circ$}}}

\newcommand{\comp}{\mathbin{\mathchoice
		{\compcent\scriptstyle}{\compcent\scriptstyle}
		{\compcent\scriptscriptstyle}{\compcent\scriptscriptstyle}}}

\newcommand{\SL}{\mathrm{SL}}
\newcommand{\GL}{\mathrm{GL}}

\newcommand{\cO}{\mathfrak{o}}
\newcommand{\mat}[4]{\left[ \begin{matrix}
#1 & #2 \\  #3 & #4 \end{matrix}\right] }

\newcommand{\Z}{\mathbb Z}

\newcommand{\val}{\mathrm{val}}

\renewcommand{\wp}{\mathfrak{p}}
\renewcommand{\det}{{\mathrm{det}}}

\newcommand{\Char}{\mathrm{Char}}

\newcommand{\bol}[1]{\bold{#1}}
\newcommand{\G}{\mathrm{G}}

\newcommand{\IR}{\mathrm{IR}}

\newcommand{\pr}{\mathrm{pr}}
\newcommand{\calp}{\mathcal{P}}
\newcommand{\dvr}{\mathrm{DVR}}
\newcommand{\dvrtwo}{\dvr_2}
\newcommand{\dvrp}{\dvr_p}
\newcommand{\dvrtwoplus}{{\dvrtwo^+}}
\newcommand{\dvrtwozero}{{\dvrtwo^\circ}}
\newcommand{\ee}{{\mathrm{e}}}

\title[Construction of Representations of $\SL_2(\cO)$]{Construction of Representations of Compact Special \\ Linear Groups of degree two}

\author{M Hassain}
\address{ MH:   The Institute of Mathematical Sciences (HBNI), CIT Campus, Taramani, Chennai 600113, India.}

\email{hassainm@hri.res.in}



\begin{document}
\begin{abstract}
	We construct the finite-dimensional continuous complex representations of $\SL_2$ over compact discrete valuation rings of even residual characteristic. We also prove that the complex  group algebras of $\SL_2$ over finite quotient rings of such compact discrete valuation rings depend on the characteristic of the ring. In particular, we prove that  the group algebras  $\mathbb C[\SL_2( \Z/2^{r}\Z)]$ and  $\mathbb C[\SL_2(\mathbb F_2[t]/(t^{r}))]$ are not isomorphic for any $r \geq 4.$
	
	%We construct the finite-dimensional continuous complex representations of $\SL_2$ over the ring of integers of non-Archimedean local fields of even residual characteristic. We also prove that the group algebras  $\mathbb C[\SL_2( \Z/2^{r}\Z)]$ and  $\mathbb C[\SL_2(\mathbb F_2[t]/(t^{r}))]$ are not isomorphic for any $r \geq 4.$\hfoot{change abstract}

\end{abstract}
\maketitle

\tableofcontents

\section{Introduction}

Let $\cO$ be a compact discrete valuation ring with maximal ideal $\wp$ and finite residue field $\mathbb F_q$ of characteristic $p.$ For $\ell\geq 1,$ the finite quotient $\cO/\wp^\ell$ is denoted by $\cO_\ell.$
For a prime $p$, let $\dvrp$ denote the set of all compact discrete valuation rings with finite residue field  of characteristic $p.$ Let  $\dvrp^{\circ} = \{ \cO \in \dvrp \mid \Char(\cO) = 0   \}$ and  $\dvrp^{+} = \{ \cO \in \dvrp \mid \Char(\cO) = p   \}$. For $\cO \in \dvrp^{\circ}$, we say $\cO$ has ramification index $\ee$ if 
$p \cO = \pi^\ee \cO,$ where $\pi$ is a fixed uniformizer of the ring $\cO,$ that is $ \wp=\pi \cO.$ 

 Let $\GL_n(\cO) $ be the group of $n \times n$ invertible matrices with entries from $\cO.$ Let $\SL_n(\cO)$ be the subgroup of $\GL_n(\cO)$ consisting of all determinant one matrices.
For $\G=\GL_n \text{ or } \SL_n,$  it is well known that for every finite-dimensional continues complex irreducible representation $\rho$ of $\G(\cO),$ there exists a smallest natural number $r=r(\rho)$ such that $\rho$ factor through the principal congruence quotient $\G(\cO_{r+1}).$ In this case, we say $\rho$ is a representation of level $r.$

In this article we focus on the construction of %level $r$
the finite-dimensional continuous complex 
irreducible representations of groups $\SL_2(\cO).$ %for every $r\geq1.$ 
For $\cO \in \dvrp$ with $p\neq2,$ the construction  is already known, see~Jaikin-Zapirain~\cite{MR2169043}. 
For $\cO=\mathbb Z_2,$ a construction of irreducible representations of $\SL_2(\cO)$ has appeared in a series of articles by Nobs and Nobs-Wolfart \cite{MR0444787,  MR0429742,  MR0387432, MR0444788}. 
Recently, in \cite{mypaper1}, we  constructed all 
 $(2r-1)\text{-level}$ irreducible representations of  $\SL_2( \cO)$ for all  $\cO \in \dvrtwoplus$ with $r \geq 1$, and for all $\cO \in \dvrtwozero$ with ramification index $\ee$ and $r \geq 2\ee.$
We also studied the representation growth of $\SL_2( \cO)$ and  proved that for $\cO \in \dvrtwoplus,$ the abscissa of convergence of the representation zeta function of  $\SL_2( \cO)$ is $1,$ resolving the last remaining open
case of this problem. By using the construction, we additionally proved that for $\cO \in \dvrtwozero$ with ramification index $\ee$ and $\cO' \in \dvrtwoplus,$   contrary
to the expectation, the group algebras  $\mathbb C[\SL_2(\cO_{2r})]$ and  $\mathbb C[\SL_2( \cO'_{2r})]$
are not isomorphic for any $r > \ee$.

In this article we give a construction of all 
%primitive 
$2r\text{-level}$ irreducible representations of  groups
$\SL_2( \cO)$ for all  $\cO \in \dvrtwoplus$ with $r \geq 1$, and for $\cO \in \dvrtwozero$ with ramification index $\ee$ and $r \geq 2\ee$
(see Theorems~\ref{thm:construction-SL-r odd- char=0 and 2(half)} and \ref{thm:construction-SL-r odd- char=2 (non inv trace)}
). This completes
the construction for all 
%finite-dimensional continues complex 
irreducible representations of $\SL_2(\cO)$ for $\cO \in \dvrtwo.$  
  We use Clifford Theory (see Isaacs~\cite[Theorem~6.11]{MR2270898}) as our main tool for the construction.
  % of the primitive irreducible representations of $\SL_2(\cO_r)$.
 From the construction we obtain information regarding the dimensions of irreducible representations. This  helps us to prove the following result, which is an extended version of \cite[Theorem 1.2]{mypaper1}.
Recall that if  $\cO \in \dvrtwoplus$ and $\cO /\wp \cong \mathbb F_q,$ then $\cO \cong  \mathbb F_q[[t]].$  
 
\begin{thm}
	\label{thm:group-algebras}
	Let $\cO \in \dvrtwozero$ with ramification index $\ee$, $\cO' \in \dvrtwoplus$ such that $\cO/\wp \cong \cO'/\wp'.$  Then the group algebras $\mathbb C[\SL_2(\cO_{r})]$ and  $\mathbb C[\SL_2( \cO'_{r})]$ are not isomorphic for any $r \geq 2 \mathrm{e}+2.$ 
\end{thm}

 In particular: 
\begin{corollary}
	The group algebras  $\mathbb C[\SL_2(\mathbb Z/(2^{r} \mathbb Z))]$ and  $\mathbb C[\SL_2(\mathbb F_2[t]/(t^{r}))]$ are not isomorphic for any $r \geq 4.$  
\end{corollary}

    See Section~\ref{Sec:group-algebras}  for its proof. 
     Since this article is a continuation of  \cite{mypaper1}, 
     %by M. Hassain and P. Singla, 
    we assume all the definitions and notations from  \cite{mypaper1} here. 
 
 \subsection*{Acknowledgement:}  The author is greatly thankful  to  Pooja Singla  for her exceptional support, without that this paper and the research behind it would not have been possible.

\section{Basic framework and outline of ideas} 
\label{sec: basic-framework} 
The basic ideas discussed in this section are already appeared  in \cite{mypaper1}. We include those here for the completeness of this article.
Let  $\cO \in \dvrp$ and $r \geq 1.$ For any $i \leq r,$ let  $f_{r,i}: M_2(\cO_{r})  \rightarrow M_2(\cO_i)$ be the natural surjective ring homomorphism.  The restriction of $f_{r,i}$ to $\GL_2(\cO_r),$ denoted by $\rho_{r,i},$ defines a surjective group homomorphism from $\GL_2(\cO_r)$ onto $\GL_2(\cO_i)$ and $\rho'_{r,i}: \SL_2(\cO_{r})  \rightarrow \SL_2(\cO_i)$ is the corresponding homomorphism obtained by restricting $f_{r,i}$ to $\SL_2(\cO_{r}) .$ For any $A \in M_2(\cO_{r}),$ the image $f_{r,1}(A)$ is denoted by $\bar{A}.$ Let $M^i =  \mathrm{ker} (\rho_{r,i})$ and $K^i = \mathrm{ker}(\rho'_{r,i}),$ then it is clear that $K^i = M^i \cap \SL_2(\cO_r).$ These are called the $i^{th}$ congruence subgroups of $\GL_2(\cO_r)$ and $\SL_2(\cO_r)$ respectively.

 For $\G=\GL_2 \text{ or } \SL_2,$ 
 an irreducible representation $\rho$ of $\G(\cO_r)$ is called a {\it twist} of another representation $\rho'$  if there exists a one-dimensional representation $\chi$ of $\G(\cO_r)$ such that $\rho' \cong \rho\otimes \chi.$ An irreducible representation $\rho$ of $\G(\cO_r)$ is called {\it primitive} if neither $\rho$ nor any of its twists 
 are trivial when restricted to $(r-1)^{th}$ congruence subgroup and $\rho$  is called {\it imprimitive} if it is not primitive. The set of imprimitive representations of $\G(\cO_r)$ can be constructed from the representations of  $\G(\cO_{r-1}).$ 
 So to understand irreducible representations of $\G(\cO),$ it is enough to understand the primitive irreducible representations of $\G(\cO_{r}),$ for $r\geq 1.$ 
 %Recall the well known result that for every finite-dimensional continues complex irreducible representation $\phi$ of $\G(\cO),$ there exists a positive natural number $r$ such that $\phi$ factor through the principal congruence quotient $\G(\cO_{r}).$

 The group of one-dimensional representations of an abelian group $H$ is denoted by $\widehat{H}.$ 
For $r \geq 2,$ fix $\ell = \lceil r/2 \rceil$ and $\ell' = \lfloor r/2 \rfloor.$ 
Let $\pi$ be a fixed uniformizer of the ring $\cO.$ We fix an additive  one-dimensional representation $\psi: \cO_r  \rightarrow \mathbb C^\times$ such that $\psi(\pi^{r-1}) \neq 1.$ 
For  $A \in M_2(\cO_{\ell'}),$ let $\tilde{A} \in M_2(\cO_r)$ be an arbitrary lift of $A$ satisfying $f_{r,\ell'}(\tilde{A}) = A.$ Define $\psi_A: %I + \pi^{r-i} M_2(\cO_r)
M^\ell \rightarrow \mathbb C^\times$ by $\psi_A(I + \pi^{\ell} B) = \psi(\pi^{\ell}\mathrm{trace}(\tilde{A}B))$ for all $I + \pi^{\ell} B \in M^{\ell}.$  Then $\psi_A$ is a well defined one-dimensional representation of $M^{\ell}=I + \pi^{\ell} M_2(\cO_r).$ 
Further  we obtain  
\begin{equation*}
%\label{eq:duality}
M_2(\cO_{\ell'}) \cong \widehat {M^{\ell}}\,;\,\,\,\, A \mapsto \psi_A.
\end{equation*}
Since $M^\ell$ is abelian, the set $\{\psi_A|_{K^\ell} \mid A \in M_2(\cO_{\ell'}) \}$ forms the set of all irreducible representations of the subgroup $K^\ell.$
It is easy to see that $\psi_A|_{K^\ell} = \psi_B|_{K^\ell}$ if and only if $A = B +xI$ for some $x \in \cO_{\ell'}.$  Let $\sim$ be  the later equivalence relation on $M_2(\cO_{\ell'})$
% as  $A \sim B$ if and only if $A = B +xI$ for some $x \in \cO_{\ell'}$
 and the equivalence class of $A \in M_2(\cO_{\ell'})$ is denoted by $[A].$ 
 Then $[A] \mapsto \psi_A|_{K^\ell}$ gives a bijection between $M_2(\cO_{\ell'})/\sim$ and $\widehat{K^\ell}.$ In view of this, we denote $\psi_A|_{K^\ell}$ by $\psi_{[A]}.$

Recall that a matrix $A \in M_2(\cO_r)$ for $r \geq 1,$ is called cyclic if there exists a vector $v \in \cO_r^{\oplus n}$  such that $\{v, Av, A^{2}v, \ldots, A^{n-1} v\}$ generate $\cO_r^{\oplus n}$ as a  free $\cO_r$-module. 
As an application of Nakyama's lemma we see that a matrix $A \in M_2(\cO_\ell)$ is cyclic if and only if $\bar{A}$ is.
% Motivated by this, we say one-dimensional representations $\psi_A$ and $\psi_{[A]}$ of groups $M^\ell$ and $K^\ell$ respectively are cyclic if and only if $A$ is cyclic.
 For  an irreducible representation $\rho$ of $\GL_2(\cO_r),$ it is easy to observe that   $\rho$ is primitive if and only if $\langle \rho|_{M^\ell}, \psi_{A} \rangle \neq 0 $ implies $A$ is a cyclic matrix. A similar criteria also holds for 
 irreducible representations of $\SL_2(\cO_r)$ as well.  

%A representation $\phi: G \rightarrow \GL(V)$ of group $G$ is denoted by $(\phi, V).$ Whenever there is no ambiguity we also denote $(\phi,V)$ by either $\phi$ or $V$ itself. 
For a subgroup $H$ of $G$ and $h \in H,$ let $C_G(h) = \{g \in G \mid gh = hg   \}$ denote the centralizer of $h$ in $G.$ Similarly for any representation $\phi$ of a normal subgroup $N$ of $G,$ the group $C_G(\phi) = \{ g \in G \mid \phi^g \cong \phi \}$ denote the inertia group (also called stabilizer) of $\phi.$ For a group $G$ and an irreducible representation $\phi$ of a subgroup $H$ of $G,$ let $\mathrm{Irr}(G \mid \phi)$ denote the set of all inequivalent irreducible constituents of the induced representation $\mathrm{Ind}_H^G(\phi).$ 
%We loosely also call $\mathrm{Irr}(G \mid \phi)$ to be the set of irreducible representations of $G$ lying above $\phi.$ 

By Clifford theory, %(see \cite[Theorem~6.11]{MR2270898}),
 for every $\psi_A \in \widehat{M^\ell}$ and $\psi_{[A]} \in \widehat{K^\ell},$ the following sets are in bijection (via induction). 
\begin{equation}
\label{eq:gl-cliiford}
\mathrm{Irr}( C_{\GL_2(\cO_{r})}(\psi_{A}) \mid \psi_A)   \leftrightarrow \mathrm{Irr}( \GL_2(\cO_r) \mid \psi_A) 
\end{equation}
\begin{equation}
\label{eq:sl-clifford}
 \mathrm{Irr}(C_{\SL_2(\cO_{r})}(\psi_{[A]}) \mid \psi_{[A]})  \leftrightarrow \mathrm{Irr}(\SL_2(\cO_{r}) \mid \psi_{[A]}).  
\end{equation} 
Therefore to understand $\mathrm{Irr}(\SL_2(\cO_{r}) \mid \psi_{[A]}),$ we will concentrate on $\mathrm{Irr}(C_{\SL_2(\cO_{r})}(\psi_{[A]}) \mid \psi_{[A]}).$

Let $A\in M_2(\cO_{\ell'})$ be  cyclic. As far as we concerned about $\mathrm{Irr}(\SL_2(\cO_{r}) \mid \psi_{[A]}),$ we can assume  $A = \mat 0 {a^{-1}\alpha } a\beta $ for some $\alpha ,\beta \in \cO_{ \ell'}$ and $a\in \cO_{ \ell'}^\times$ (see \cite[Section~4]{mypaper1}).
Fix a lift $\tilde{A} =  \mat 0 {\tilde{a}^{-1} \tilde{\alpha}}{\tilde a}{\tilde{\beta}}  \in M_2(\cO_{r})$ of  $A.$
%$ = \mat 0 {a^{-1}\alpha } a\beta \in M_2(\cO_{\ell'})$.
 For $i \in \{ \ell, \ell' \}$, define 
\[ \mathrm{h}_{\tilde{A}}^{i} = \{ x \in \cO_{r } \mid 2x = 0  \,\, \mathrm{mod}\, (\pi^{i} ), \, \,  x(x+\tilde \beta) = 0 \,\, \mathrm{mod}\, (\pi^{i} ) \}.
\]
Let  $\mathrm{H}_{\tilde{A}}^{i} = \{ e_x = \mat1{\tilde{a}^{-1}x}01 \mid x \in h_{\tilde{A}}^i \}  $. Then $\mathrm{H}_{\tilde{A}}^{i} $ is an abelian group for $i \in \{ \ell, \ell' \}$.
 We remark that for the definition of $h_{\tilde{A}}^i;$ $i=\ell' , \ell$ by a lift of $A=\mat 0 {a^{-1}\alpha } a\beta$ we always consider a matrix of the form $\mat 0 {\tilde{a}^{-1} \tilde{\alpha}}{\tilde a}{\tilde{\beta}}$.
 % It suffices for our purposes to consider lifts of this kind, see Remark~\ref{remark:lift-dependence}. See Section~\ref{sec:Stabilizer-of-Sl2-action} for a few results regarding  $h_{\tilde{A}}^{\ell'}$ and $h_{\tilde{A}}^{\ell}$. 
Recall the following lemma from \cite{mypaper1}.

\begin{lemma}
	\label{lem:stabilizer-form} \cite[Lemma 2.2]{mypaper1}
	Let $\cO \in \dvrp$ and  $A = \mat 0 {a^{-1}\alpha } a\beta \in M_2(\cO_{\ell'})$ be cyclic. Then  
	\begin{enumerate}	
		\item  $ C_{\GL_2(\cO_{r})} (\psi_A)  = C_{\GL_2(\cO_{r})} (\tilde{A}) M^{\ell'}.$ 
		\item For $p \neq 2,$  $ C_{\SL_2(\cO_{r})} (\psi_{[A]})  = C_{\SL_2(\cO_{r})}(\psi_A).$ 
		\item For $p=2,$ $C_{\SL_2(\cO_{r})} (\psi_{[A]})  = C_{\SL_2(\cO_{r})} (\psi_{A}) \mathrm{H}_{\tilde{A}}^{\ell'}.$ 	
	\end{enumerate}	
\end{lemma}
The group $C_{\GL_2(\cO_{r})} (\tilde{A})$ is abelian for cyclic $A$. Hence the one-dimensional representations $\psi_A$ and $\psi_{[A]}$ extend to the groups $C_{\GL_2(\cO_{r})} (\tilde{A}) M^{\ell}$ and $C_{\GL_2(\cO_{r})} (\tilde{A}) M^{\ell} \cap \SL_2(\cO_r)$ respectively (see \cite[Lemmas~3.1 and 3.2]{mypaper1}). Now onwards, as in \cite{mypaper1}, we use the following notations for simplification:

\begin{itemize} 
	\item $C_G^\ell(\tilde{A}) :=  C_{\GL_2(\cO_{r})} (\tilde{A}) M^{\ell}$ and $C_G^{\ell'}(\tilde{A}) := C_{\GL_2(\cO_{r})} (\tilde{A}) M^{\ell'}.$
	\item  $C_S^\ell(\tilde A):= C_G^\ell(\tilde{A}) \cap \SL_2(\cO_{r})$ and $C_S^{\ell'}(\tilde A) := C_G^{\ell'}(\tilde{A}) \cap \SL_2(\cO_{r}).$
\end{itemize} 
 We also use the following new notations.
\begin{itemize} 
	\item $D_G^\ell(\tilde{A}):=(C_{\GL_2(\cO_{r})} (\tilde{A})\cap M^1) M^{\ell}$ and   $D_G^{\ell'}(\tilde{A}):=(C_{\GL_2(\cO_{r})} (\tilde{A})\cap M^1) M^{\ell'}$
	\item $D_S^{\ell}(\tilde A):=D_G^{\ell}(\tilde{A})\cap \SL_2(\cO_{r})$ and $D_S^{\ell'}(\tilde A):=D_G^{\ell'}(\tilde{A})\cap \SL_2(\cO_{r})$
\end{itemize}
Note that $D_G^{\ell}(\tilde{A}),$  $D_G^{\ell'}(\tilde{A}),$  $D_S^{\ell}(\tilde A)$ and  $D_S^{\ell'}(\tilde A) $ are the first congruent subgroups of $C_G^{\ell}(\tilde{A}),$ $C_G^{\ell'}(\tilde{A})$ $C_S^{\ell}(\tilde{A})$ and $C_S^{\ell'}(\tilde{A})$ respectively.
Also $ K^{\ell}\trianglelefteq D_S^{\ell}(\tilde A) \trianglelefteq D_S^{\ell'}(\tilde A)\trianglelefteq C_S^{\ell'}(\tilde A)\trianglelefteq C_{\SL_2(\cO_{r})}(\psi_{[A]}).$

Let  $\cO \in \dvrtwo.$ We construct the primitive irreducible representations of $\SL_2(\cO_{r})$ for odd $r.$  Note that for odd $r,$ we have $r=2\ell'+1.$ 
 For the construction we consider the following cases separately.
\begin{enumerate}
	\item $\cO \in \dvrtwozero$ and $\ell'\geq 2\mathrm{e}+1.$ 
	\item $\cO \in \dvrtwoplus$ and $\ell'\geq 1$ with $\mathrm{trace}(A)\in \cO_{\ell'}^\times.$
	\item $\cO \in \dvrtwoplus$ and $\ell'\geq 1$ with $\mathrm{trace}(A)\in \pi\cO_{\ell'}.$
\end{enumerate}

For the first and second cases, the construction  
% of    $\mathrm{Irr}( C_{\SL_2(\cO_{2\ell'+1})} (\psi_{[A]})\mid \psi_{[A]} )$ 
is given by the following theorem.
%
%\begin{thm}\label{thm:construction-SL-r odd- char=0}
%	Let $\cO \in \dvrtwozero$
%	and $\ell' \geq 2\mathrm{e}+1.$ Let $A \in M_2(\cO_{\ell'})$ be cyclic. For every $\delta\in \mathrm{Irr}( C_{\SL_2(\cO_{2\ell'+1})} (\psi_{[A]})\mid \psi_{[A]}),$ there exists a one-dimensional representation  $\phi \in \mathrm{Irr}(D_S^\ell(\tilde A)\mid \psi_{[A]}) $ and a $q$-dimensional representation  $\rho_{\phi} \in \mathrm{Irr}(C_S^{\ell'}(\tilde A)\mid \phi)$ such that $\delta \cong \mathrm{Ind}_{C_S^{\ell'}(\tilde A)}^{ C_{\SL_2(\cO_{2\ell'+1})} (\psi_{[A]})}(\rho_\phi).$ 
%
%\end{thm}

\begin{thm}\label{thm:construction-SL-r odd- char=0 and 2(half)}
	Let $\cO \in \dvrtwozero$
	and $\ell' \geq 2\mathrm{e}+1$ (resp.  $\cO \in \dvrtwoplus$
	and $\ell' \geq 1.$)  The following hold for every cyclic $A \in M_2(\cO_{\ell'})$  (resp. every cyclic $A \in M_2(\cO_{\ell'})$ such that $\mathrm{trace}(A)\in \cO_{\ell'}^\times$).

\begin{enumerate}
	\item Every representation $\phi \in \mathrm{Irr}(D_S^\ell(\tilde A)\mid \psi_{[A]}) $  is one-dimensional.
	\item For each $\phi \in \mathrm{Irr}(D_S^\ell(\tilde A)\mid \psi_{[A]}), $ 	there exists exactly one representation, say $\hat{\phi},$ in $\mathrm{Irr}(D_S^{\ell'}(\tilde A)\mid  \phi), $ and $\dim (\hat{\phi})=q.$
	\item For $\phi \in \mathrm{Irr}(D_S^\ell(\tilde A)\mid \psi_{[A]}), $ the representation $\hat{\phi}$ extends to $C_S^{\ell'}(\tilde A).$ Further each $\rho \in \mathrm{Irr}(C_S^{\ell'}(\tilde A)\mid \hat{\phi})$ is an extension of $\hat{\phi}.$ %has dimension $q.$
	\item By Clifford theory, the set of equivalence classes of $\mathrm{Irr}(C_S^{\ell'}(\tilde A)\mid  \psi_{[A]})$ under the conjugation action of $C_{\SL_2(\cO_{2\ell'+1})} (\psi_{[A]})$ is in bijective correspondence with $\mathrm{Irr}(\SL_2(\cO_{2 \ell'+1})\mid  \psi_{[A]})$ given by  $\rho \mapsto \mathrm{Ind}_{C_S^{\ell'}(\tilde A)}^{ \SL_2(\cO_{2\ell'+1})}(\rho).$
	
\end{enumerate}
	
\end{thm}

%
%
%For the second case, the construction of is given by the following theorem.
%\begin{thm}\label{thm:construction-SL-r odd- char=2}
%	Let $\cO \in \dvrtwoplus$
%	and $\ell' \geq 1.$ Let $A \in M_2(\cO_{\ell'})$ be cyclic such that $\mathrm{trace}(A)\in \cO_{\ell'}^\times.$ For every $\delta\in \mathrm{Irr}( C_{\SL_2(\cO_{2\ell'+1})} (\psi_{[A]})\mid \psi_{[A]}),$ there exists a one-dimensional representation  $\phi \in \mathrm{Irr}(D_S^\ell(\tilde A)\mid \psi_{[A]}) $ and a $q$-dimensional representation  $\rho_{\phi} \in \mathrm{Irr}(C_S^{\ell'}(\tilde A)\mid \phi)$ such that $\delta \cong \mathrm{Ind}_{C_S^{\ell'}(\tilde A)}^{ C_{\SL_2(\cO_{2\ell'+1})} (\psi_{[A]})}(\rho_\phi).$ 
%	%Further $\mathrm{Ind}_{C_S^{\ell'}(\tilde A)}^{ C_{\SL_2(\cO_{2\ell'+1})} (\psi_{[A]})}(\rho_\phi)\cong \mathrm{Ind}_{C_S^{\ell'}(\tilde A)}^{ C_{\SL_2(\cO_{2\ell'+1})} (\psi_{[A]})}(\rho'_{\phi'})$ if and only if $\phi \cong \phi'$ and $\rho_{\phi} \cong \rho'_{\phi'}.$
%	
%	
%	
%\end{thm}

For the third case, the construction is given by the following theorem.

\begin{thm}\label{thm:construction-SL-r odd- char=2 (non inv trace)}
	Let  $\cO \in \dvrtwoplus$
	and $\ell' \geq 1.$  The following hold for  every cyclic $A \in M_2(\cO_{\ell'})$ such that $\mathrm{trace}(A)\in \pi \cO_{\ell'}.$

	\begin{enumerate}
		\item Every representation $\phi \in \mathrm{Irr}(D_S^\ell(\tilde A)\mid \psi_{[A]}) $  is one-dimensional.
		\item For each $\phi \in \mathrm{Irr}(D_S^\ell(\tilde A)\mid \psi_{[A]}), $ 	there exists a group $\mathbb M_\phi$ such that 
		\begin{enumerate}
			\item $D_S^\ell(\tilde A) \leq \mathbb M_\phi \leq  C_{\SL_2(\cO_{2\ell'+1})} (\psi_{[A]}) .$ 
			%and $[M_\phi: D_S^\ell(\tilde A)]\geq q.$
			\item $\phi$ extends to $ \mathbb M_\phi.$  Further each $\rho \in \mathrm{Irr}(\mathbb M_\phi\mid \phi)$ is an extension of $\phi.$
			\item  For each $\rho \in \mathrm{Irr}(\mathbb M_\phi\mid \phi),$ the induced representation $\mathrm{Ind}_{\mathbb M_\phi}^{ \SL_2(\cO_{2\ell'+1})}(\rho)$ is irreducible.

		\end{enumerate}

	\end{enumerate}
	
\end{thm}

We prove Theorem~\ref{thm:construction-SL-r odd- char=0 and 2(half)}  in Section~\ref{sec:proof of construction cases 1&2}
and Theorem~\ref{thm:construction-SL-r odd- char=2 (non inv trace)} in Section~\ref{sec:proof of construction cases 3}. Note that  as far as we concerned about $\mathrm{Irr}(\SL_2(\cO_{r}) \mid \psi_{[A]}),$ to prove Theorems~\ref{thm:construction-SL-r odd- char=0 and 2(half)} 
and \ref{thm:construction-SL-r odd- char=2 (non inv trace)},  we can assume  $A = \mat 0 {a^{-1}\alpha } a\beta $ for some $\alpha ,\beta \in \cO_{ \ell'}$ and $a\in \cO_{ \ell'}^\times$ (see \cite[Section~4]{mypaper1}).
% For $\beta \in \pi\cO_{ \ell'}$ we can further assume $\alpha \in \pi\cO_{ \ell'}$ (see \cite[Proposition~4.2]{mypaper1}).
We also fix a lift $\tilde{A} =  \mat 0 {\tilde{a}^{-1} \tilde{\alpha}}{\tilde a}{\tilde{\beta}}  \in M_2(\cO_{r})$ of  $A.$

In Section~\ref{sec:characterization-of-E*a-results}, we define the subset  $\mathbb E^\prime_{\tilde{A}}$ of $\mathrm{H}_{\tilde{A}}^{\ell'}$ and give a characterization of it. We see that the set  $\mathbb E^\prime_{\tilde{A}}$ plays a very important role in proving Theorems~\ref{thm:construction-SL-r odd- char=0 and 2(half)}  and \ref{thm:construction-SL-r odd- char=2 (non inv trace)}.

%%%%%%%%%%%%%%%%%%%%%%%%%%%%%%%%%%%
\section{Definition and Characterization of elements of $\mathbb E^\prime_{\tilde{A}}$} 
\label{sec:characterization-of-E*a-results}

%Let $\cO \in \dvrtwo. $
In this section, we define the subset  $\mathbb E^\prime_{\tilde{A}}$ of $\mathrm{H}_{\tilde{A}}^{\ell'}$ and give a characterization of it. The set  $\mathbb E^\prime_{\tilde{A}}$ plays a very important role in the construction of primitive irreducible representations of $\SL_2(\cO_{r})$  for odd $r.$ 
 Let $\cO \in \dvrtwo. $ Throughout this section, we fix $A = \mat 0{a^{-1} \alpha}a{\beta} \in M_2(\cO_{\ell'})$ and its lift $\tilde{A} = \mat 0{\tilde{a}^{-1}\tilde \alpha}{\tilde{a}}{\tilde \beta} \in M_2(\cO_r).$ For $\lambda \in \cO_r$, the matrix $\mat{1}{\tilde{a}^{-1} \lambda}{0}{1}\in M_2(\cO_r)$ is denoted by $e_\lambda.$ 
%The set  $\mathbb E^\prime_{\tilde{A}}$ plays a very important role in the construction of primitive irreducible representations of $\SL_2(\cO_{r})$ for odd $r.$
 Recall that $\mathrm{H}_{\tilde{A}}^i = \{ e_\lambda  \mid 2\lambda = 0  \,\, \mathrm{mod}\, (\pi^{i} ), \, \,  \lambda(\lambda+\tilde \beta) = 0 \,\, \mathrm{mod}\, (\pi^{i} ) \},$ for $i\in\{\ell, \ell'\}.$

From \cite[Lemma~3.2]{mypaper1}, we obtain that $\psi_{[A]}$ extends to $C_S^\ell(\tilde{A})$ and more over every representation in  $ \mathrm{Irr}(C_S^{\ell}(\tilde A)\mid \psi_{[A]})$ is an extension of $\psi_{[A]}.$
Since $D_S^\ell(\tilde{A})\leq C_S^\ell(\tilde{A}),$ similar result also holds for $D_S^\ell(\tilde{A}).$
As in \cite{mypaper1}, define 
\[
\mathbb E_{\tilde{A}} := \{e_\lambda  = \mat 1{\tilde{a}^{-1} \lambda} 01 \in \mathrm{H}_{\tilde{A}}^\ell \mid \psi_{[A]}\,\, \mathrm{extends}\,\, \mathrm{to}\,\, C_S^{\ell}(\tilde{A}) \langle e_\lambda\rangle\},
\]
where $\langle e_\lambda\rangle$ denotes the group generated by $ e_\lambda.$ 
Now we define the new set 
\[
\mathbb E^\prime_{\tilde{A}} := \{e_\lambda  = \mat 1{\tilde{a}^{-1} \lambda} 01 \in \mathrm{H}_{\tilde{A}}^{\ell'} \mid \psi_{[A]}\,\, \mathrm{extends}\,\, \mathrm{to}\,\, D_S^{\ell}(\tilde{A}) \langle e_\lambda\rangle\}.
\]
Since $\mathrm{H}_{\tilde{A}}^{\ell}\subseteq \mathrm{H}_{\tilde{A}}^{\ell'}$ and $D_S^{\ell}(\tilde{A})\leq C_S^{\ell}(\tilde{A}),$ we have $	\mathbb{E}_{\tilde{A}}\subseteq 	\mathbb{E}^\prime_{\tilde{A}}.$ More precisely,  $	\mathbb{E}_{\tilde{A}}=	\mathbb{E}^\prime_{\tilde{A}} \cap \mathrm{H}_{\tilde{A}}^{\ell}$ (see Corollary~\ref{cor:EtildeA and E'tildeA}).
%The set  $\mathbb E^\prime_{\tilde{A}}$ plays a very important role in the construction of primitive irreducible representations of $\SL_2(\cO_{r})$ for odd $r.$  Note that for odd $r,$ we have $r=2\ell'+1.$ 
A characterization of the elements in $\mathbb E_{\tilde{A}}$ is given in \cite{mypaper1}. By modifying those results, we characterize the elements of $\mathbb E^\prime_{\tilde{A}}$ here.
%
%
%
%
%
%in Section~\ref{sec:characterization-of-E*a-results}.
%We use the characterization results of $\mathbb E^\prime_{\tilde{A}}$  to give a construction of all primitive irreducible representations of $\SL_2(\cO_{2 \ell'+1})$. 
%
%In this section, we give a characterization of $\mathbb E^\prime_{\tilde{A}}.$ 
%Throughout this section, we fix $A = \mat 0{a^{-1} \alpha}a{\beta} \in M_2(\cO_{\ell'})$ and its lift $\tilde{A} = \mat 0{\tilde{a}^{-1}\tilde \alpha}{\tilde{a}}{\tilde \beta} \in M_2(\cO_r).$ For $\lambda \in \cO_r$, the matrix $\mat{1}{\tilde{a}^{-1} \lambda}{0}{1}\in M_2(\cO_r)$ is denoted by $e_\lambda.$ 
%Recall that  $\mathbb E_{\tilde{A}} = \{e_\lambda   \in \mathrm{H}_{\tilde{A}}^\ell \mid \psi_{[A]}\,\, \mathrm{extends}\,\, \mathrm{to}\,\, C_S^{\ell}(\tilde{A})\langle e_\lambda\rangle\}$ and  $\mathbb E^\prime_{\tilde{A}} = \{e_\lambda   \in \mathrm{H}_{\tilde{A}}^{\ell'} \mid \psi_{[A]}\,\, \mathrm{extends}\,\, \mathrm{to}\,\, D_S^{\ell}(\tilde{A})\langle e_\lambda\rangle\}.$
%
%
%%
%
%
%
%
The sets  
$\mathbb E_{\tilde{A}}$ and $\mathbb E^\prime_{\tilde{A}}$ are in bijective correspondence with the sets $ E_{\tilde{A}} := \{\lambda \in \mathrm{h}_{\tilde{A}}^\ell \mid e_\lambda \in \mathbb E_{\tilde{A}} \}$ and $ E^\prime_{\tilde{A}} := \{\lambda \in \mathrm{h}_{\tilde{A}}^{\ell'} \mid e_\lambda \in \mathbb E^\prime_{\tilde{A}} \}$ respectively. Recall that $\mathrm{h}_{\tilde{A}}^i = \{ x \in \cO_{r } \mid 2x = 0  \,\, \mathrm{mod}\, (\pi^{i} ), \, \,  x(x+\tilde \beta) = 0 \,\, \mathrm{mod}\, (\pi^{i} ) \},$ for $i\in\{\ell, \ell'\}.$

For $x, y, \lambda \in \cO_r,$ define
\begin{eqnarray*}
	f(\lambda, x, y) &=&xy \lambda (\tilde \beta- \lambda ) -\tilde \alpha \lambda y^2 + \lambda ( x^2-1) ,\\
	g(x, y) &=& x^2 + \tilde{\beta} xy - \tilde{\alpha} y^2.
\end{eqnarray*}
Note that $	f(\lambda, x, y)=\lambda ( g(x,y)-1) - \lambda^2  xy .$
We will keep these notations fixed throughout this article.

For $\lambda \in \mathrm{h}_{\tilde{A}}^{\ell'} $
and fixed lifts $\tilde{\alpha}, \tilde{\beta} \in \cO_r$ of $\alpha, \beta \in \cO_{\ell'},$ we define the set  
\begin{equation*} 
%\label{eq: *} 
E'_{\lambda, \tilde A} = \{ (x, y) \in \cO_r \times \pi\cO_r \mid g(x,y) = 1 \,\, \mathrm{mod}\,(\pi^\ell) \,\, \mathrm{and} \,\, \lambda y \in \pi^\ell \cO_r  \}
\end{equation*} 
and its subset $E_{\lambda, \tilde A}^{\prime\circ} = \{ (x, y) \in E_{\lambda, \tilde{A}}  \mid  \psi(f(\lambda, x, y)) = 1 \}.  
$ Note that the sets $E'_{\lambda, \tilde A}$ and $E_{\lambda, \tilde A}^{\prime\circ}$ are related to the corresponding sets $E_{\lambda, \tilde A}$ and $E_{\lambda, \tilde A}^{\circ}$ in \cite{mypaper1} by $E'_{\lambda, \tilde A}=E_{\lambda, \tilde A}\cap (\cO_r \times \pi\cO_r )$ and $E_{\lambda, \tilde A}^{\prime\circ}=E^{\circ}_{\lambda, \tilde A}\cap (\cO_r \times \pi\cO_r ).$
The following result explores elements of $E^\prime_{ \tilde A}$ through various equivalent conditions. It is a modified version of \cite[Theorem~5.2]{mypaper1}.  
\begin{thm}
	\label{thm:condition for extension} 
	For $\lambda \in \mathrm{h}_{\tilde{A}}^{\ell'},$ the following are equivalent. 
	\begin{enumerate}
		\item $E^{\prime}_{\lambda, \tilde A} = E_{\lambda, \tilde A}^{\prime\circ} .$
		\item Any element of $K^\ell$ which is of the form $[e_\lambda, X]$ for $X \in D_{S}^\ell(\tilde{A})$ is in the kernel of $\psi_{[A]}.$ 
		\item $[\langle e_\lambda\rangle, D_{S}^\ell(\tilde{A}) ] \cap K^\ell$ is contained in the kernel of $\psi_{[A]}.$ 
		\item  $ \lambda \in E^\prime_{ \tilde A},$ i.e. there exists an extension of $\psi_{[A]} $ to $D_{S}^\ell(\tilde{A})\langle e_\lambda\rangle.$
	\end{enumerate}
\end{thm}
To prove Theorem~\ref{thm:condition for extension}, we need the following lemma. We also use this result in Section~\ref{sec:proof of construction cases 1&2}.

\begin{lemma}\label{lem:odd r- construction-GL}
	Let  $A \in M_2(\cO_{\ell'})$ be cyclic and $\tilde{A}   \in M_2(\cO_{r})$  be a lift  of $A.$ 
	Then  the following hold.
	\begin{enumerate}
		
		\item The one-dimensional representation $\psi_{A}$ of $M^{\ell}$  extends to $D_G^{\ell}(\tilde A),$   
		and
		every representation $\phi \in  \mathrm{Irr}(D_G^{\ell}(\tilde A)\mid  \psi_{A}) $ is one-dimensional. 
		\item Each $\phi \in  \mathrm{Irr}(D_G^{\ell}(\tilde A)\mid  \psi_{A}) $ is stabilised by $C_G^{\ell'}(\tilde A).$ 
	
	\end{enumerate}
\end{lemma} 

\begin{proof} 
	First of all note that the subgroups $D_G^{\ell}(\tilde A),$ $D_G^{\ell'}(\tilde A)$ and $C_G^{\ell'}(\tilde A)$ are independent of the lift $\tilde{A}.$
	Since $D_G^{\ell}(\tilde A)\leq C_G^{\ell}(\tilde A),$ (1) directly follows from \cite[Lemma~3.1(4)]{mypaper1}. 
	Let $\phi \in  \mathrm{Irr}(D_G^{\ell}(\tilde A)\mid  \psi_{A}) .$ Since $C_G^{\ell'}(\tilde A)=C_{\GL_2(\cO_{r})} (\tilde{A})M^{\ell'},$  to show (2) it is enough to show that $\phi$ is stabilised by both $C_{\GL_2(\cO_{r})} (\tilde{A})$ and  $M^{\ell'}.$ For that let $ I+\pi X \in C_{\GL_2(\cO_{r})} (\tilde{A})\cap M^1, $ $  I+\pi^{\ell}B\in M^{\ell},$ $Y \in C_{\GL_2(\cO_{r})} (\tilde{A})$ and $ I+\pi^{\ell'}D\in M^{\ell'}.$
	Since $C_{\GL_2(\cO_{r})} (\tilde{A})$ abelian and $UV=VU$ for all $U\in M^{\ell}$ and $V\in M^{\ell'},$ we have
	$$ [Y,( I+\pi X )( I+\pi^{\ell}B)]=( I+\pi X )[Y, I+\pi^{\ell}B]( I+\pi X )^{-1}\in M^{\ell}\,\, \mathrm{and}$$
	$$[ I+\pi^{\ell'}D,( I+\pi X )( I+\pi^{\ell}B)] =[ I+\pi^{\ell'}D, I+\pi X ]= I+\pi^{\ell'}(D-( I+\pi X )D( I+\pi X )^{-1})\in M^{\ell'+1}\subseteq M^{\ell}.$$
	Therefore 
	\begin{eqnarray*}
	\phi([Y,( I+\pi X )( I+\pi^{\ell}B)])&=&\psi_A(( I+\pi X )[Y, I+\pi^{\ell}B]( I+\pi X )^{-1})\\
	&=&\psi_A([Y, I+\pi^{\ell}B])\\
	&=&1
	\end{eqnarray*}
%$\phi([Y,( I+\pi X )( I+\pi^{\ell'+1}B)])=\psi_A(( I+\pi X )[Y, I+\pi^{\ell'+1}B]( I+\pi X )^{-1})=1$
 and
\begin{eqnarray*}
\phi([ I+\pi^{\ell'}D,( I+\pi X )( I+\pi^{\ell'+1}B)]) &=&\psi_A( I+\pi^{\ell'}(D-( I+\pi X )D( I+\pi X )^{-1}))\\
&=&\psi(\pi^{\ell'}\mathrm{trace}(\tilde A(D-( I+\pi X )D( I+\pi X )^{-1})))\\
&=&\psi(\pi^{\ell'}\mathrm{trace}(\tilde AD)-\pi^{\ell'} \mathrm{trace}( \tilde A( I+\pi X )D( I+\pi X )^{-1}))\\
&=&\psi(\pi^{\ell'}\mathrm{trace}(\tilde AD)-\pi^{\ell'} \mathrm{trace}(( I+\pi X )\tilde AD( I+\pi X )^{-1}))\\
&=&\psi(\pi^{\ell'}\mathrm{trace}(\tilde AD)-\pi^{\ell'} \mathrm{trace}(\tilde AD))\\
&=&\psi(0)=1.
\end{eqnarray*}
	 % $\phi([ I+\pi^{\ell'}D,( I+\pi X )( I+\pi^{\ell'+1}B)]) =\psi_A( I+\pi^{\ell'}(D-( I+\pi X )D( I+\pi X )^{-1}))=1.$ 
	 Hence $\phi$ is stabilised by both $C_{\GL_2(\cO_{r})} (\tilde{A})$ and  $M^{\ell'}.$ 
	
\end{proof}

\begin{proof}[Proof of Theorem~\ref{thm:condition for extension} ]
		 This theorem can be proved by a simple modification in the proof of \cite[Theorem~5.2]{mypaper1}. We first  replace the groups $C_{S}^\ell(\tilde{A}) $ and $C_{G}^{\ell}(\tilde{A}) $  by their subgroups $D_{S}^\ell(\tilde{A}) $ and $D_{G}^{\ell}(\tilde{A}) $ respectively. For example, to prove (2) implies (3), we  consider an extension $\widehat{\psi_{A}}$  of $\psi_{[A]}$ to the group $D_{G}^\ell(\tilde{A}),$ not to the group $C_{G}^\ell(\tilde{A}).$ We also replace $ \mathrm{h}_{\tilde{A}}^{\ell}$ by $ \mathrm{h}_{\tilde{A}}^{\ell'}.$ It is easy to complete the proof by  using  the following observations. So we skip the complete proof here.
		\begin{itemize}
			\item Every $ X\in D_{G}^\ell(\tilde{A})$ is of the form $X = ((1+\pi x)I + \pi y \tilde A)Z$ for some $x,y \in \cO_r$ and  $Z \in M^\ell .$ 
			\item For $x,y \in \cO_r,$ observe that  $g(1+\pi x,\pi y)=1 \,\,\mathrm{mod}\, (\pi^\ell)$  if and only if there exists $X \in  D_{S}^\ell(\tilde{A})$ such that $X = ((1+\pi x)I + \pi y \tilde A)Z$ for some $Z \in M^\ell .$
			\item For  $ \lambda \in \mathrm{h}_{\tilde{A}}^{\ell'},$ from \cite[Lemma~5.1(1)]{mypaper1} we have $e_\lambda \tilde{A} e_\lambda^{-1}= \tilde{A} + \lambda I \,\, \mathrm{mod}\,(\pi^{\ell'}).$ Hence  $ e_\lambda  \pi\tilde{A} e_\lambda^{-1}= \pi \tilde{A} + \pi\lambda I \,\, \mathrm{mod}\,(\pi^{\ell}).$ This gives $ e_\lambda ^c X e_\lambda^{-c}  \in D_{G}^\ell(\tilde{A}) ,$ for all $X\in D_{G}^\ell(\tilde{A})$ and $c\geq 0.$
			\item For  $ \lambda \in \mathrm{h}_{\tilde{A}}^{\ell'},$ $e_\lambda^2 \in M^{\ell'} \subseteq C_{G}^{\ell'}(\tilde{A}).$ Hence by Lemma~\ref{lem:odd r- construction-GL}(2), $e_\lambda^{2}$ stabilizes  $\widehat{\psi_{A}}.$  Therefore, for every $X\in D_{G}^\ell(\tilde{A})$ and $c\geq 0,$ 
			$$\widehat{\psi_{A}} (e_\lambda^{c} X e_{\lambda}^{-c} )=\begin{cases} \widehat{\psi_{A}} (e_\lambda X e_{\lambda}^{-1} ) &  \mathrm{if}\,\, c \,\, \mathrm{is}\,\, \mathrm{odd},\\ \widehat{\psi_{A}} ( X ) &  \mathrm{if}\,\, c \,\, \mathrm{is}\,\, \mathrm{even}.
			\end{cases}$$
		\end{itemize}

\end{proof}
The following corollary gives the relation between the sets $E_{\tilde{A}}$ and $E^\prime_{\tilde{A}}.$
\begin{corollary}\label{cor:EtildeA and E'tildeA}
	$E_{\tilde{A}}=E^\prime_{\tilde{A}} \cap \mathrm{h}_{\tilde{A}}^{\ell}.$ In particular, for even $r,$ $E_{\tilde{A}}=E^\prime_{\tilde{A}}.$
\end{corollary}
\begin{proof}
	By definition of $E_{\tilde{A}}$ and $	E^\prime_{\tilde{A}} ,$ it is easy to observe that $E_{\tilde{A}}\subseteq	E^\prime_{\tilde{A}} .$ Therefore we must have 	$	E_{\tilde{A}}\subseteq 	E^\prime_{\tilde{A}} \cap \mathrm{h}_{\tilde{A}}^{\ell}.$ To show 	$E^\prime_{\tilde{A}} \cap \mathrm{h}_{\tilde{A}}^{\ell} \subseteq 	E_{\tilde{A}},$ let $\lambda \in E^\prime_{\tilde{A}} \cap \mathrm{h}_{\tilde{A}}^{\ell}.$  If $\lambda \in \pi^{\ell}\cO_r,$ then $C_S^{\ell}(\tilde{A}) \langle e_\lambda\rangle=C_S^{\ell}(\tilde{A}).$ Since $C_S^{\ell}(\tilde{A})\leq C_G^{\ell}(\tilde{A}),$  by \cite[Lemma~3.1(4)]{mypaper1}, there exists an extension of $ \psi_{[A]}$ to $C_S^{\ell}(\tilde{A}).$ Therefore we obtain $\lambda \in 	E_{\tilde{A}}.$
So assume $\lambda \notin \pi^{\ell}\cO_r.$   To show $\lambda \in 	E_{\tilde{A}},$ by  \cite[Theorem~5.2]{mypaper1} it is enough to show the following.
\begin{equation}\label{eqn:csl-dsl}
\{[e_\lambda, X]\mid X\in C_S^{\ell}(\tilde{A})\}\cap K^\ell \subseteq \mathrm{ker}(\psi_{[A]}).
\end{equation}
%$\{[e_\lambda, X]\mid X\in C_S^{\ell}(\tilde{A})\}\cap K^\ell \subseteq \mathrm{ker}(\psi_{[A]}).$ 
We claim that $\{[e_\lambda, X]\mid X\in C_S^{\ell}(\tilde{A})\}\cap K^\ell = \{[e_\lambda, X]\mid X\in D_S^{\ell}(\tilde{A})\}\cap K^\ell.$ Since $\lambda \in E^\prime_{\tilde{A}} ,$ the claim along with Theorem~\ref{thm:condition for extension} gives (\ref{eqn:csl-dsl}).

To show the claim,
%let  $\lambda \in \mathrm{h}_{\tilde{A}}^{\ell}\setminus(\pi^{\ell})$ and   
$X \in C_{S}^\ell(\tilde{A})$ be such that   $[e_\lambda, X]\in K^\ell.$ By definition of $C_{S}^\ell(\tilde{A}),$ there exists $x I+y\tilde A \in C_{\GL_2(\cO_{r})} ( \tilde A)$ and $ I+\pi^\ell B\in M^\ell$ such that $X=(x I+y\tilde A)( I+\pi^\ell B).$ So we have 
$[e_\lambda, x I+y\tilde A]=[e_\lambda, X]= I\,\, \mathrm{mod}\, (\pi^\ell).$ Therefore by \cite[Lemma~5.1(2)]{mypaper1}, we  have $\lambda y \in \pi^{\ell}\cO_r.$ Since $\lambda \notin \pi^{\ell}\cO_r,$ we obtain $y\in \pi\cO_r.$
Combining this with $\det(X)=1$  (because  $X %=(x I+y\tilde A)( I+\pi^\ell B)
\in C_{S}^\ell(\tilde{A})$), we get 
$$ 1= \det(X)=x^2\,\, \mathrm{mod}\, (\pi).$$
Therefore $x=1\,\, \mathrm{mod}\, (\pi).$ Hence $X =(x I+y\tilde A)( I+\pi^\ell B)\in D_{S}^\ell(\tilde{A}).$ This gives the claim.
	
\end{proof}
\begin{corollary}\label{cor:E_tilde A remark}
	The following hold for $\mathrm{h}_{\tilde{A}}^{\ell'}$ and   $E^\prime_{\tilde{A}}.$ 
	\begin{enumerate}
		\item 	 $  \pi^{\ell'}\cO_r\subseteq E^\prime_{\tilde{A}} \subseteq \mathrm{h}_{\tilde{A}}^{\ell'} .$ 	
		
		\item If $\lambda \in E^\prime_{\tilde{A}},$ then $\lambda+\pi^{\ell'}\cO_r\subseteq E^\prime_{\tilde{A}} .$ 
	\end{enumerate}
\end{corollary}

\begin{proof}
	For  $\lambda \in \pi^{\ell'}\cO_r,$ we have $E^\prime_{\lambda, \tilde{A}} = E_{\lambda, \tilde{A}} ^{\prime\circ}.$ This follows because for $(x,y) \in E^\prime_{\lambda, \tilde{A}} ,$ we must have $g(x,y) = 1 \,\, \mathrm{mod}\,(\pi^\ell)$ and $ \lambda y \in \pi^{\ell}\cO_r$. Therefore, 
	$	f(\lambda,x,y)	
	= \lambda (  g(x,y)-1) - \lambda ^2 xy  = 0.$
%	\begin{eqnarray*}
%		f(\lambda,x,y)	
%		&=&xy \lambda (\tilde \beta- \lambda ) -\tilde \alpha \lambda y^2 + \lambda ( x^2-1)\\
%		& =& \lambda (  g(x,y)-1) - \lambda ^2 xy \\
%		&=&  \lambda ( g(x,y)-1) - \lambda( \lambda y) x  = 0.
%	\end{eqnarray*}
	So (1) follows by Theorem~\ref{thm:condition for extension} . 
	To show (2), let $\lambda \in E^\prime_{\tilde{A}}$ and $z\in  \pi^{\ell'}\cO_r.$ 
	By definition of $\mathrm{h}_{\tilde{A}}^{\ell'},$  it is easy to observe that $\lambda+z \in \mathrm{h}_{\tilde{A}}^{\ell'}.$ Since $\lambda \in E^\prime_{\tilde{A}},$  by Theorem~\ref{thm:condition for extension} , we have 
	$E^\prime_{\lambda, \tilde{A}}=E^{\prime\circ}_{\lambda, \tilde{A}}.$ Therefore if we show $E^\prime_{\lambda+z, \tilde{A}}=E^\prime_{\lambda, \tilde{A}}$ and $E^{\prime\circ}_{\lambda+z, \tilde{A}}=E^{\prime\circ}_{\lambda, \tilde{A}},$ then we have 
	$\lambda+z \in E^\prime_{\tilde{A}}. $ 
	
	Note that  $zy\in \pi^{\ell}\cO_r$ for all $y \in \pi\cO_r.$ So by definition of $E^\prime_{\lambda, \tilde{A}},$ it is easy to show that 
	$E^\prime_{\lambda+z, \tilde{A}}=E^\prime_{\lambda, \tilde{A}}.$ Next for $(x,y) \in E^\prime_{\lambda+z, \tilde{A}} ,$ we have 
	\begin{eqnarray*}
		f(\lambda + z,x,y)	&=& (\lambda + z)  ( g(x,y)-1 ) - (\lambda + z)^2 xy\\
		&=& \lambda   ( g(x,y)-1 ) - \lambda^2  xy\\
		&=&f(\lambda ,x,y).
	\end{eqnarray*}
	Here the second equality follows because $z ( g(x,y)-1 )=z\lambda y=z^2 y=0.$ %\hfoot{is it clear?}
	Therefore by definition of $E^{\prime\circ}_{\lambda, \tilde{A}},$ we have  $E^{\prime\circ}_{\lambda+z, \tilde{A}}=E^{\prime\circ}_{\lambda, \tilde{A}}.$
	Hence (2) holds.
\end{proof}
%
%.............good remark...but not usefull in this article.....
%
%\begin{remark}\label{rmk:cor:E_tilde A remark}
%For $\lambda \in \mathrm{h}_{\tilde{A}}^{\ell'}$ and $z\in  \pi^{\ell'}\cO_r,$ we have $\lambda+z \in \mathrm{h}_{\tilde{A}}^{\ell'}$ and 
%	$E^\prime_{\lambda+z, \tilde{A}}=E^{\prime}_{\lambda, \tilde{A}}.$ Further $	f(\lambda + z,x,y)=f(\lambda ,x,y)$ for all $(x,y)\in E^{\prime}_{\lambda, \tilde{A}}.$
%\end{remark}
%
%
%...................................
%%
Let $\boldsymbol{\psi}:\mathbb F_q \rightarrow \mathbb C^\times$ be an additive one-dimensional representation given by
$\boldsymbol{\psi}(\bar x)=\psi(\pi^{r-1} x),$  for all $x\in \cO_r$ such that $x=\bar x \,\, \mathrm{mod}\,(\pi).$ By \cite[Lemma~5.4]{mypaper1}, there exists a unique $\xi \in \mathbb F_q^\times$ such that $\mathrm{ker}(\boldsymbol{\psi})=\{ \xi x^2 +x \mid x \in \mathbb F_q\}.$

\begin{definition}\label{defn:of xi}
	({\bf Definition of $\xi$}) Now onwards, for  fixed $\psi : \cO_r \rightarrow \mathbb C^\times$ and $\boldsymbol{\psi}:\mathbb F_q \rightarrow \mathbb C^\times$ as above, we fix the notation $\xi$ for the unique element in $\mathbb F_q^\times$ such that
	$\mathrm{ker}(\boldsymbol{\psi})=\{ \xi x^2 +x \mid x \in \mathbb F_q\}.$
\end{definition}
Let $\cO \in \dvrtwozero. $
Let $\ee$ be the ramification index of $\cO.$ The characterization of  $E^\prime_{\tilde{A}}$ for $r > 4\ee$ is given by the following result.  It is analogous to \cite[Theorem~5.6]{mypaper1}.

\begin{thm}\label{S_A-in-number}
	Let $\cO \in  \dvrtwozero$. The following hold for $E^\prime_{\tilde{A}}$. 
	\begin{enumerate}
		\item  For $ r \geq 2(\ee+1)$ with $\beta \in \cO_{\ell'}^\times,$ 	$E_{\tilde{A}}^{\prime} =\pi^{\ell'} \cO_r .$	
		\item  For $ r > 4\ee$ with $\beta\in \pi \cO_{\ell'},$	$E_{\tilde{A}}^{\prime} =\pi^{\ell'}\cO_r .$

	\end{enumerate}
\end{thm}
\begin{proof}%\hfoot{add proper ref and ideas}
	
	\noindent {\bf (1): Let $r \geq 2(\ee+1)$ and $\beta \in \cO_{\ell'}^\times.$} 
	Note that $\ell' > \ee.$	By  \cite[Proposition~3.8(2)]{mypaper1}, we have
	$	\mathrm{h}_{\tilde{A}}^{\ell'}= \pi^{\ell'}\cO_r .$ Therefore by  Corollary~\ref{cor:E_tilde A remark}(1), we obtain that
	$E_{\tilde{A}}^\prime=\pi^{\ell'}\cO_r .$
	
	\noindent {\bf (2): Let $r > 4 \ee$ and $\beta \in \pi \cO_{\ell'}.$}
	From \cite[Proposition~3.8(2)]{mypaper1}, we have 
	$
	\mathrm{h}_{\tilde{A}}^{\ell'}=\pi^{\ell'-m_0}\cO_r ,$ where $ m_0= \min \{\ee, \val(\tilde\beta)\}.$ By Corollary~\ref{cor:E_tilde A remark}(1), we  have $\pi^{\ell'}\cO_r \subseteq E_{\tilde{A}}^\prime.$ Therefore by Theorem~\ref{thm:condition for extension}, to show $E_{\tilde{A}}^\prime=\pi^{\ell'}\cO_r ,$ it is enough to prove that for any $\lambda \in 	\mathrm{h}_{\tilde{A}}^{\ell'} \setminus \pi^{\ell'}\cO_r  ,$ there exists 
	$(x,y ) \in E_{\lambda, \tilde A}^\prime\setminus E_{\lambda, \tilde A}^{\prime\circ}.$
	Let  $\lambda \in 	\mathrm{h}_{\tilde{A}}^{\ell'} \setminus \pi^{\ell'}\cO_r .$ Then  $\lambda = \pi^{i}u,$ for some $u \in \cO_{r}^\times$ and $\ell'-m_0\leq  i \leq \ell'-1.$
	Since $\ee$ is the ramification index of $\cO,$ there exists $w\in \cO_r^\times$ such that $2=\pi^\ee w.$
	Let  $x=1+\pi^{r-\ee-i-1}(wu)^{-1}$ and $y=0.$
Now we show that this choice of $(x,y)$ is in $E_{\lambda, \tilde A}^\prime\setminus E_{\lambda, \tilde A}^{\prime\circ}. $
	
	\noindent {\bf \underline{$\lambda y = 0 \,\, \mathrm{mod}\,(\pi^\ell) $:}} This is obvious. 
	
	\noindent {\bf \underline{$g(x, y) = 1 \,\, \mathrm{mod}\,(\pi^ \ell)$: }}  For this we note that 
	$$g(x, y) =1 +2\pi^{r-\ee-i-1}(wu)^{-1}+\pi^{2r-2\ee-2i-2}(wu)^{-2}=1 +\pi^{r-i-1}u^{-1}+\pi^{2r-2\ee-2i-2}(wu)^{-2}.$$
	So it remains to prove $r-i-1\geq \ell$ and $2r-2\ee-2i-2\geq \ell. $  By using the facts that $i\leq \ell' - 1 $ and $\ell\geq 2\ee+1,$   we obtain
	$$r-i-1\geq (\ell+\ell')-(\ell'-1)-1 =\ell \, \, \,  \mathrm{and}$$
	$$2r-2\ee-2i-2\geq 2(\ell+\ell') - (\ell-1) -2(\ell'-1)-2=\ell+1\geq \ell. $$
	
	\noindent{\bf \underline{$\psi (f(\lambda, x, y))\neq 1$: }} By the definition of $\psi,$ it is enough to show that 
	$f(\lambda, x, y) =\pi^{r-1} .$ For the given values of $x$ and $y,$ we have 
	\begin{eqnarray*}
		f(\lambda, x, y) &=&\lambda(x^2-1)\\
		&=&\pi^{i}u( 2\pi^{r-\ee-i-1}(wu)^{-1}+\pi^{2r-2\ee-2i-2}(wu)^{-2})\\
		&=& \pi^{r-1}+\pi^{2r-2\ee-i-2}w^{-2}u^{-1}.
	\end{eqnarray*}
Further, by using the facts that $ \ell\geq 2\ee+1$ and $i\leq \ell' - 1,$   we obtain
$$2r-2\ee-i-2\geq 2r-(\ell-1)-(\ell'-1)-2=2r-(\ell+\ell')=r.$$
Therefore $	f(\lambda, x, y) =\pi^{r-1}.$
	Hence  $ (x,y) \in  E_{\lambda, \tilde A}^\prime\setminus E_{\lambda, \tilde A}^{\prime\circ}.$ This completes the proof.
	% of Theorem~\ref{S_A-in-number}.

\end{proof}
	
We now proceed to 
characterize the elements of $E^\prime_{\tilde{A}}$ for $\cO \in \dvrtwoplus.$  For that we fix a few notations.  These notations will be used throughout  the article for $\cO \in \dvrtwoplus$. We assume $\cO \in \dvrtwoplus$ for rest of this section.
%\hfoot{check the validity}

For $v \in \cO_r,$  and $0\leq i \leq  r-1,$ we denote $(v)_i\ $ for the unique elements in $\mathbb F_q $ such that $v= (v)_0 + (v)_1 \pi +\cdots + 
(v)_{r-1} \pi^{r-1}.$ Note that  $v \in \cO_r^2$
if and only if $(v)_{2i+1} = 0$ for all $i$.  

\begin{definition} ({\bf Definition of $\psi$ for $\cO \in \dvrtwoplus$})
	Let $\boldsymbol{\psi}:\mathbb F_q \rightarrow   \mathbb C^\times$ be a fixed additive one-dimensional representation such that $\boldsymbol{\psi}(1)\neq 1.$  Define $\psi: \cO_r  \rightarrow   \mathbb C^\times$ by  $\psi(v)=\boldsymbol{\psi}((v)_{r-1})$ for  $v\in \cO_r.$
\end{definition}

For $A = \mat 0 {a^{-1} \alpha} a \beta \in M_2(\cO_{\ell'}),$ we define (as in \cite{mypaper1})  $\bol k$ and $\bol s$ as follows.   
\begin{eqnarray}
\bol  k & = & \val(\beta). \nonumber \\
\bol s & = &\begin{cases}
2 \lfloor \bol k/2 \rfloor + 1 &  \mathrm{if}\,\,   \alpha = v^2 \,\,\mathrm{mod}\,(\pi^ \bol k ), \\ m &  \mathrm{if}\,\, \alpha = v_1^2 + \pi^m v_2^2 \,\,\mathrm{mod}\,(\pi^ \bol k ) \,\,\mathrm{for}\,\,\mathrm{odd}\,\, m < \bol k \,\,\mathrm{and}\,\, v_2 \in \cO_{\ell'}^\times. 
\end{cases}\nonumber
\end{eqnarray}
We remark that for any lift $\tilde{ \beta}\in\cO_r$ of $\beta,$ we have $\bol k=  \min \{ \val({\tilde{\beta}}), \ell' \}.$ 
%The following examples illustrate these notions. 
Recall the following lemma from \cite{mypaper1}.
\begin{lemma}\cite[Lamma~5.12]{mypaper1}
	\label{defintion-of-s-and-w1} 
	For any lift $\tilde{\alpha}$ of $\alpha,$ there exists  $w_1, w_2 \in \cO_{r}$ such that $\tilde \alpha = w_1^2 + \pi^{\bol s} w_2^2 \in \cO_r.$ Further, if $\bol s < \bol k$ then $w_2 \in \cO_{r}^\times$. 
\end{lemma}

\begin{remark}\label{proper-tilde-alpha-and-tilde-beta}
	({\bf Definition of $\tilde{\alpha},$ $\tilde{\beta},$ $w_1$ and $w_2.$}) Now onwards, we fix lifts $\tilde{\alpha}$ and $\tilde{\beta}$ of $\alpha$ and $\beta$ respectively. We also fix $w_1, w_2 \in \cO_r$  such that  $\tilde{\alpha}= w_1^2 + \pi^{\bol s} w_2^2 \in \cO_r.$ For $A = \mat 0{a^{-1} \alpha}a{\beta}$, we fix $\tilde{A}  = \mat 0{\tilde{a}^{-1}\tilde \alpha}{\tilde{a}}{\tilde \beta} $, where $\tilde \alpha$ and  $\tilde \beta$ correspond to the above fixed choices. 
\end{remark}
For $\lambda \in \cO_r$ and $\tilde{A}  = \mat 0{\tilde{a}^{-1}\tilde \alpha}{\tilde{a}}{\tilde \beta} , $ we define $\bol i_{\lambda, \tilde{A}}$, $\bol j_{\lambda, \tilde{A}}$ and $\delta_{\lambda, \tilde{A}}$.  
\begin{eqnarray}
\bol i_{\lambda, \tilde{A}}& = &\val({\lambda}), \nonumber \\
\bol j_{\lambda, \tilde{A}} & = & \min \{ \val({\lambda+ \tilde{\beta}}), \ell' \}, \nonumber \\
\delta_{\lambda, \tilde{A}}&=& \bol j_\lambda - \bol s - \mathrm{max} \{ \ell - \bol i_\lambda, \ell - \bol k, \lceil (\ell - \bol s) /2 \rceil \}. \nonumber
\end{eqnarray}
We remark that $\tilde{A}  = \mat 0{\tilde{a}^{-1}\tilde \alpha}{\tilde{a}}{\tilde \beta} $ is fixed throughout, therefore we leave  $\tilde{A}$ notation from  $\bol i_{\lambda, \tilde{A}}, \bol j_{\lambda, \tilde{A}}, \delta_{\lambda, \tilde{A}}$ and write these as $\bol i_\lambda$, $\bol j_{\lambda}$, $\delta_\lambda$ respectively. We also define $\epsilon$ as follows:
$$\epsilon=\begin{cases}
1& \mathrm{if}\, r=\mathrm{even}\\0 & \mathrm{if} \, r=\mathrm{odd.}
\end{cases}$$
 Note that these notations are same as in \cite{mypaper1}. We will use these  throughout the article.

Let $\cO_{r}^2 =\{ x^2 \mid x \in \cO_r\}.$  Note that from Corollary~\ref{cor:E_tilde A remark}(1), we have $\pi^{\ell'}\cO_r \subseteq E^\prime_{\tilde{A}}.$ 
For $\lambda \in  \mathrm{h}_{\tilde{A}}^{\ell'} \setminus \pi^{\ell'}\cO_r ,$ the following  result, which is  analogous to \cite[Theorem~5.18]{mypaper1},  discusses the conditions on $\lambda$ such that $ \lambda \in E^\prime_{\tilde{A}}.$ 
% We include\hfoot{To do} a proof of this result in Appendix~\ref{appendix}. 

\begin{thm} 
	\label{thm:simplification of extension conditions} 
	
	For $\lambda \in  \mathrm{h}_{\tilde{A}}^{\ell'} \setminus \pi^{\ell'}\cO_r ,$ $\lambda \in E_{\tilde{A}}^\prime $ if and only if the following conditions hold:
	\begin{enumerate}[(I)]
		\item $ \lambda \in \pi^{\ell - \ell'}\cO_{r}^2 \,\, \mathrm{mod}\,  (\pi^{\ell'}).$
		\item$2\bol j_\lambda + \bol i_\lambda = 2 \ell' +\bol  s - \epsilon .$ 
		\item For  $ \bol j_\lambda < \ell',$ $\bol s < \bol k$ and $\delta_\lambda \geq 0,$
		we must have 
		$
		\xi u_1^2 u_2^2 =  u_1 w_2^2 \,\, \mathrm{mod}\, (\pi^{2 \delta_\lambda +1}), 
		$ 
		where $u_1,u_2 \in \cO_{r}^\times$ such that $\lambda = \pi^{\bol i_\lambda} u_1$
		and $ \lambda + \tilde{\beta} = \pi^{\bol j_\lambda} u_2.$
	\end{enumerate}
\end{thm} 
To prove Theorem~\ref{thm:simplification of extension conditions}, we need the following lemma.
\begin{lemma}\label{lem:link lemma}
		For $\lambda \in  \mathrm{h}_{\tilde{A}}^{\ell'} \setminus \pi^{\ell'}\cO_r ,$ $\lambda \in E_{\tilde{A}}^\prime $ implies $\lambda \in  \mathrm{h}_{\tilde{A}}^{\ell}. $
\end{lemma}
\begin{proof}
	 For even $r,$ since $\ell =\ell',$ the lemma holds trivially. So assume $r$ is odd. Let
	 $\lambda \in  \mathrm{h}_{\tilde{A}}^{\ell'} \setminus \pi^{\ell'}\cO_r $ be such that  $\lambda \in E_{\tilde{A}}^\prime .$ We first claim that $\lambda \in \pi \cO_{r}.$ If $  \lambda \in  \cO_{r}^\times,$ then choose $\mu \in \cO_r^\times$ such that $\mu^2 = \lambda^{-1}\,\, \mathrm{mod}\, (\pi).$ Existence of such $\mu \in \cO_r^\times$ is because the map $\bar{z}\mapsto \bar{z}^2$ is a bijection from 
	 $\cO_1$ to itself. Let $x=1+\pi^{\ell'}\mu$ and $y=0.$ Now we show that this choice of $(x,y)$ is in $E_{\lambda, \tilde A}^\prime\setminus E_{\lambda, \tilde A}^{\prime\circ}. $ 
	 
	 \noindent {\bf \underline{$\lambda y = 0 \,\, \mathrm{mod}\,(\pi^\ell) $:}} This is obvious. 
	 
	 \noindent {\bf \underline{$g(x, y) = 1 \,\, \mathrm{mod}\,(\pi^ \ell)$: }}  This follows because
	 $g(x, y) =1 +\pi^{2\ell'}\mu^2 $ and $2\ell' \geq\ell.$ 
	
	 \noindent{\bf \underline{$\psi (f(\lambda, x, y))\neq 1$: }} By the definition of $\psi,$ it is enough to show that 
	 $f(\lambda, x, y) =\pi^{r-1} .$ For the given values of $x$ and $y,$ we have 
	 $$ f( \lambda,x,y) = \lambda \pi^{2\ell'}\mu^2 = \pi^{r-1}.$$   
	 The last equality follows because $r=2\ell'+1$ and $\pi^{2\ell'}\mu^2 =\pi^{2\ell'}\lambda^{-1} . $ Therefore  $ (x,y) \in  E_{\lambda, \tilde A}^\prime\setminus E_{\lambda, \tilde A}^{\prime\circ}.$ Hence by Theorem~\ref{thm:condition for extension}, $\lambda \notin E_{\tilde{A}}^\prime .$ It is a contradiction to our hypothesis that $\lambda \in E_{\tilde{A}}^\prime .$ Therefore  $  \lambda \notin  \cO_{r}^\times$ and hence the claim follows.
	 
Next we proceed to show $\lambda \in  \mathrm{h}_{\tilde{A}}^{\ell}. $  Since $\lambda \in  \mathrm{h}_{\tilde{A}}^{\ell'},$  by definition of $\mathrm{h}_{\tilde{A}}^{\ell'},$  we have $ \lambda(\tilde \beta + \lambda) = 0 \,\, \mathrm{mod}\, (\pi^{\ell'} ).$ 
If $\lambda \notin \mathrm{h}_{\tilde{A}}^{\ell}, $ then $ \lambda(\tilde \beta + \lambda) \neq 0 \,\, \mathrm{mod}\, (\pi^{\ell} ).$ Therefore there exists $u \in \cO_{r}^\times$ such that $\lambda(\tilde \beta + \lambda)=\pi^{\ell'}u.$ Since $\lambda\notin \pi^{\ell'}\cO_r,$ we must have $\tilde \beta + \lambda\in \pi \cO_r.$ Hence the claim $ \lambda\in \pi \cO_r$ implies $ \tilde \beta \in  \pi \cO_r.$ 
Let $y=\pi^{\ell'}u^{-1}$ and $x=1+y w_1.$ We now proceed to prove that this choice of $(x,y)$ is in $E_{\lambda, \tilde A}^\prime\setminus E_{\lambda, \tilde A}^{\prime\circ}. $
  
 \noindent {\bf \underline{$\lambda y = 0 \,\, \mathrm{mod}\,(\pi^\ell) $:}} This follows because $y=\pi^{\ell'}u^{-1}$  and $\lambda \in \pi \cO_r.$ 
 
\noindent {\bf \underline{$g(x, y) = 1 \,\, \mathrm{mod}\,(\pi^ \ell)$: }}  For this we note that 
$$g(x, y) =1 + \tilde \beta x \pi^{\ell'} u^{-1}+w_2^2 \pi^{2\ell'} u^{-2} .$$
So it remains to prove  $\tilde \beta x \pi^{\ell'} u^{-1}+w_2^2 \pi^{2\ell'} u^{-2} = 0 \,\, \mathrm{mod}\,(\pi^ \ell). $ Since $\tilde \beta \in \pi\cO_r$ and $2\ell'\geq \ell$ we are done.

\noindent{\bf \underline{$\psi (f(\lambda, x, y))\neq 1$: }} By the definition of $\psi,$ it is enough to show that 
$f(\lambda, x, y) =\pi^{r-1} .$ For the given values of $x$ and $y,$ we have 
\begin{eqnarray*}
	f(\lambda, x, y) &=&(1+w_1 y)y\lambda(\lambda+\tilde{ \beta})+\lambda \pi^\bol s w_2^2 y^2.
\end{eqnarray*}

	Since $\lambda(\lambda+\tilde{ \beta} ) = \pi^{\ell'}u$ and $y=\pi^{\ell'}u^{-1},$ we must have $(1+w_1 y)y\lambda(\lambda+\tilde{ \beta})=\pi^{2\ell'}=\pi^{r-1}.$ Further $\lambda \in \pi\cO_r$ and $y^2 \in \pi^{2\ell'}\cO_r$ implies $\lambda \pi^\bol s w_2^2 y^2=0.$ Therefore $	f(\lambda, x, y) =\pi^{r-1}.$
Hence  $ (x,y) \in  E_{\lambda, \tilde A}^\prime\setminus E_{\lambda, \tilde A}^{\prime\circ}.$ By Theorem~\ref{thm:condition for extension}, $\lambda \notin E_{\tilde{A}}^\prime .$ It is a contradiction to our hypothesis that $\lambda \in E_{\tilde{A}}^\prime .$ Therefore we must have  $  \lambda \in  \mathrm{h}_{\tilde{A}}^{\ell}.$ Hence the lemma follows.

\end{proof}
\begin{proof}[Proof of Theorem~\ref{thm:simplification of extension conditions}]
Note that Theorem~\ref{thm:simplification of extension conditions} directly follows from \cite[Theorem~5.18]{mypaper1} if we show  
$$
\{\lambda \in  \mathrm{h}_{\tilde{A}}^{\ell'} \setminus \pi^{\ell'}\cO_r  \mid \lambda \in E_{\tilde{A}}^\prime\}= \{\lambda \in  \mathrm{h}_{\tilde{A}}^{\ell} \setminus \pi^{\ell'}\cO_r  \mid \lambda \in E_{\tilde{A}}\}.$$
Since 
$E_{\tilde{A}}\subseteq E^\prime_{\tilde{A}}$ and $ \mathrm{h}_{\tilde{A}}^{\ell} \subseteq \mathrm{h}_{\tilde{A}}^{\ell'},$ we have 
$$\{\lambda \in  \mathrm{h}_{\tilde{A}}^{\ell'} \setminus \pi^{\ell'}\cO_r  \mid \lambda \in E_{\tilde{A}}^\prime\}\supseteq \{\lambda \in  \mathrm{h}_{\tilde{A}}^{\ell} \setminus \pi^{\ell'}\cO_r  \mid \lambda \in E_{\tilde{A}}\}.$$
To show the converse, let $\lambda \in  \mathrm{h}_{\tilde{A}}^{\ell'} \setminus \pi^{\ell'}\cO_r$ such that  $ \lambda \in E_{\tilde{A}}^\prime.$ By Lemma~\ref{lem:link lemma}, we have $\lambda \in \mathrm{h}_{\tilde{A}}^{\ell}.$  Since 
$E_{\tilde{A}}=E^\prime_{\tilde{A}} \cap \mathrm{h}_{\tilde{A}}^{\ell}$ (by Corollary~\ref{cor:EtildeA and E'tildeA}), we must have $ \lambda \in E_{\tilde{A}}.$ Therefore 
$$\{\lambda \in  \mathrm{h}_{\tilde{A}}^{\ell'} \setminus \pi^{\ell'}\cO_r  \mid \lambda \in E_{\tilde{A}}^\prime\}\subseteq \{\lambda \in  \mathrm{h}_{\tilde{A}}^{\ell} \setminus \pi^{\ell'}\cO_r  \mid \lambda \in E_{\tilde{A}}\}.$$

\end{proof}

The following result will be used in Section~\ref{sec:proof of construction cases 1&2} to prove Theorem~\ref{thm:construction-SL-r odd- char=0 and 2(half)}.

\begin{proposition}\label{prop:char 2 E*A} 
	%\hfoot{is it used in full?}
	Let  $\cO \in \dvrtwoplus$ and $r\geq 1$ be odd. For $\beta \in \cO_{\ell'}^\times ,$  $E^\prime_{\tilde A}=  \pi^{\ell'}\cO_r.$
	\end{proposition}

\begin{proof}	
 Since $\beta \in \cO_{\ell'}^\times,$  by definition of $\mathrm{h}_{\tilde{A}}^{\ell'},$ we have $\mathrm{h}_{\tilde{A}}^{\ell'}=\{0, \tilde \beta\}+\pi^{\ell'}\cO_r.$ Therefore to show $E^\prime_{\tilde A}=  \pi^{\ell'}\cO_r,$ it is enough to show that $\tilde \beta \notin E^\prime_{\tilde A}$ (by Corollary~\ref{cor:E_tilde A remark}). Suppose contrary that $\tilde \beta \in E^\prime_{\tilde A}.$ Then by \textit{(I)} of  Theorem~\ref{thm:simplification of extension conditions}, we have $\tilde \beta \in \pi^{\ell-\ell'}\cO_r^2.$ Note that $\ell-\ell'=1,$ because $r$ is odd. Therefore we must have $\tilde \beta \in \pi\cO_r,$ which is a contradiction to the fact that $\tilde \beta \in \cO_{\ell'}^\times.$  Hence $\tilde \beta \notin E^\prime_{\tilde A}.$
\end{proof}

%%%%%%%%%%%%%%%%%%%%%%%%%%
\section{Proof of Theorem~\ref{thm:construction-SL-r odd- char=0 and 2(half)}}\label{sec:proof of construction cases 1&2}

In this section we prove Theorem~\ref{thm:construction-SL-r odd- char=0 and 2(half)}.
% and \ref{thm:construction-SL-r odd- char=2}. 
For that we need the following machinery.
Let $G$ be a finite group. Let $N$ be a normal subgroup of $G$ such that $G/N$ is an elementary abelian $p$-group (i.e. order of $gN$ in $G/N$ is $p$ for all $g \in G\setminus N$). Then the group $G/N$ has a structure of $\mathbb F_p$ vector space.
Let $\chi:N\rightarrow \mathbb{C}^\times$ be a one-dimensional representation of $N$ such that $C_G(\chi)=G.$ Define an alternating bilinear form  $h_\chi : G/N \times G/N \rightarrow \mathbb{C}^\times $ by  $h_\chi(g_1 N, g_2 N)=\chi([g_1 , g_2] )=\chi(g_1g_2g_1^{-1}g_2^{-1}).$ By bilinearity we mean that 
\begin{eqnarray*}
	h_\chi(g_1g_2 N, g_3 N)&=&h_\chi(g_1 N, g_3 N)h_\chi(g_2 N, g_3 N)\\
	h_\chi(g_1 N, g_2g_3 N)&=&h_\chi(g_1 N, g_2 N)h_\chi(g_1 N, g_3 N),
\end{eqnarray*}
%$h_\chi(g_1g_2 N, g_3 N)=h_\chi(g_1 N, g_3 N)h_\chi(g_2 N, g_3 N)$ and $h_\chi(g_1 N, g_2g_3 N)=h_\chi(g_1 N, g_2 N)h_\chi(g_1 N, g_3 N)$
for all $g_1, g_2, g_3 \in G.$ This follows from the facts that $[g_1g_2,g_3]=g_1[g_2,g_3]g_1^{-1} [g_1,g_3] $ and $[g_1,g_2g_3]=[g_1,g_2]g_2 [g_1,g_3] g_2^{-1}.$  An easy computation shows that $h_\chi$ is well defined. Define the set 
$$\bar{R}_\chi=\{xN \in  G/N \mid h_\chi(xN,yN)=1 \, \, \mathrm{for}\, \mathrm{all}\, y \in G\}.$$
This is called the \textit{radical} of the form $h_\chi$ and we say $h_\chi$ is \textit{non-degenerate} if $\bar{R}_\chi=\{N\}.$

\begin{proposition}\label{prop:bilinear non deg.}
	Suppose $h_\chi$ is non-degenerate. Then there exists exactly one %irreducible
	 representation in  $ \mathrm{Irr}( G\mid \chi),$ 
	%of $G$ lying above $\chi,$ 
	and it has dimension $[G:N]^{\frac{1}{2}}.$
\end{proposition}
\begin{proof}
	See \cite[Proposition~8.3.3]{MR701540}.
\end{proof}	
The following corollary is a generalisation of the above lemma.
\begin{corollary}\label{cor:bilinear form}
	Assume  $h_\chi$ is not  non-degenerate. Let $R_\chi$ be the preimage of $\bar{R}_\chi$ under the canonical projection map $G\rightarrow G/N.$  Then the following holds.
	\begin{enumerate}
		\item The one-dimensional representation  $\chi$ has an extension to $R_\chi.$
		\item Every representation $\phi \in \mathrm{Irr}( R_\chi \mid \chi)$ is one-dimensional.
		\item For each  representation $\phi \in \mathrm{Irr}( R_\chi \mid \chi),$ there exists exactly one %irreducible 
		representation in  $ \mathrm{Irr}( G\mid \phi),$ 
		%of $G$ lying above $\phi,$
		 and it has dimension $[G:R_\chi]^{\frac{1}{2}}.$
	\end{enumerate}
\end{corollary}

\begin{proof}
	To our knowledge, this result was first used by Hill in  \cite{MR1334228} for construction of regular characters of $\GL_n(\cO)$ for $\cO \in\dvrp^{\circ}.$ 
	%$\Char(\cO)=0.$  
	A proof of this result is included in for example    \cite[Corollary~3.3]{MR3743488}.
\end{proof}

 Let $\cO \in \dvrtwo$ and $r\geq 1$ be odd. Then $r=2\ell'+1$ and $\ell=\ell'+1.$
 %
 %.......new discription about A.........
 Let $A= \mat 0 {a^{-1}\alpha } a \beta  \in M_2(\cO_{\ell'})$ be  cyclic. 
 Fix a lift $\tilde{A} =  \mat 0 {\tilde{a}^{-1} \tilde{\alpha}}{\tilde a}{\tilde{\beta}}  \in M_2(\cO_{2\ell'+1})$ of  $A.$
 %%%%%%%%%%%%%
 %
 %...........old discription about A.........
 %
%Let $A \in M_2(\cO_{\ell'})$ be cyclic  and $\tilde{A}   \in M_2(\cO_{2\ell'+1})$  be a lift  of $A.$   
%
%%%%%%%%%%%%%%%%%%
%
%
By Lemma~\ref{lem:odd r- construction-GL}, every  $\chi\in  \mathrm{Irr}(D_G^{\ell}(\tilde A)\mid  \psi_{A})$ is one-dimensional and  stabilised by $C_G^{\ell'}(\tilde A).$
Therefore we can consider 
the alternating bilinear form $h_\chi : D_G^{\ell'}(\tilde A)/D_G^{\ell}(\tilde A) \times D_G^{\ell'}(\tilde A)/D_G^{\ell}(\tilde A)  \rightarrow \mathbb{C}^\times $ defined by $h_\chi (xD_G^{\ell}(\tilde A), yD_G^{\ell}(\tilde A))=\chi([x,y]).$

\begin{lemma}\label{lem:non deg. of h phi} 
	For  $\chi\in  \mathrm{Irr}(D_G^{\ell}(\tilde A)\mid  \psi_{A}),$ let $h_\chi$
	% $h_\chi : (D_G^{\ell'}/D_G^{\ell}) \times (D_G^{\ell'}/D_G^{\ell})  \rightarrow \mathbb{C}^\times $ 
	be	the alternating bilinear form defined as above. Then 
	$h_\chi $ is 
	non-degenerate.
\end{lemma}
\begin{proof}
 This lemma directly follows from a general result proved by Stasinski-Stevens, see \cite[Lemma~4.5]{MR3743488}. For reader’s convenience, we give a proof here. 
	Let $xD_G^{\ell}(\tilde A)\in D_G^{\ell'}(\tilde A)/D_G^{\ell}(\tilde A)$ such that $xD_G^{\ell}(\tilde A)\neq D_G^{\ell}(\tilde A).$ Note that we can choose $x\in M^{\ell'}.$ Let $x= I+\pi^{\ell'}B.$ Since $xD_G^{\ell}(\tilde A)\neq D_G^{\ell}(\tilde A),$ we have 	$x= I+\pi^{\ell'}B\notin D_G^{\ell}(\tilde A).$ So by definition of $D_G^{\ell}(\tilde A),$ we must have 
	$ I+\pi^{\ell'}B\notin C_{\GL_2(\cO_{2\ell'+1})} (\tilde{A})\cap M^1   \,\, \mathrm{mod}\, (\pi^{\ell} ).$ Therefore  $B\notin C_{\GL_2(\cO_{2\ell'+1})} (\tilde{A})\,\, \mathrm{mod}\, (\pi ).$
Hence by definition of $C_{\GL_2(\cO_{2\ell'+1})} (\tilde{A}),$  we have $\tilde{A}B-B\tilde{A}\neq 0 \,\, \mathrm{mod}\, (\pi ).$ Now choose $Z\in M_2(\cO_{2\ell'+1})$ be such that $\mathrm{trace}((\tilde{A}B-B\tilde{A})Z) =1 \,\, \mathrm{mod}\, (\pi ).$
Then for $y= I+\pi^{\ell'}Z\in M^{\ell'},$ we obtain
\begin{eqnarray*}
	h_\chi(xD_G^{\ell}(\tilde A), y D_G^{\ell}(\tilde A))&=&\chi([x,y])\\
	&=&\chi( I+\pi^{2\ell'}(BZ-Z B))\\
	&=&\psi(\pi^{2\ell'}\mathrm{trace}(\tilde{A}BZ-\tilde{A}ZB))\\
	&=&\psi(\pi^{2\ell'})\neq 1.
\end{eqnarray*}
Here the last equality follows from the fact that  $\mathrm{trace}(\tilde{A}BZ-\tilde{A}ZB)=\mathrm{trace}(\tilde{A}BZ-B\tilde{A}Z) =1 \,\, \mathrm{mod}\, (\pi ).$
Therefore $h_\chi$ is non-degenerate.

\end{proof}

% For  cyclic $A \in M_2(\cO_{\ell'}),$
The following lemma will discuss some properties of representations in $ \mathrm{Irr}(D_S^\ell(\tilde A)\mid \psi_{[A]}).$
We also use this result in section \ref{sec:proof of construction cases 3}.

\begin{lemma}\label{lem:D_S L phi} 
	Let  $A \in M_2(\cO_{\ell'})$ be cyclic and $\tilde{A}   \in M_2(\cO_{2\ell' +1})$  be a lift  of $A.$ Then for $\phi \in \mathrm{Irr}(D_S^\ell(\tilde A)\mid \psi_{[A]}) $ the following hold.
	\begin{enumerate}
		\item  $\phi  $ is one-dimensional.
		\item $C_G^{\ell'}(\tilde A)$ stabilizes $\phi .$ 
		\item $\phi$  extends to $D_G^\ell(\tilde A)$ and every such extension is stabilised by $C_G^{\ell'}(\tilde A).$
	\end{enumerate}
\end{lemma}
\begin{proof}
	Note that (1) follows from  \cite[Lemma 3.2]{mypaper1} and the fact that   $D_S^\ell(\tilde A)\leq C_S^\ell(\tilde A).$ It is easy to observe that (3) implies (2).
	%if we prove (3), then (2) follows trivially.
	% Since $D_S^{\ell'} (\tilde A)\leq C_G^{\ell'}(\tilde A),$ it is easy to see that (3) implies (2). 
	To show (3), note that $\mathrm{Irr}(M^\ell\mid  \psi_{[A]} )=\{\psi_{A+x I} \mid x\in \cO_{\ell'} \}.$ Therefore 
	\begin{equation}\label{eqn:abcds}
\mathrm{Irr}(D_G^\ell(\tilde A)\mid  \psi_{[A]} )=\bigcup_{x\in \cO_{\ell'}}\mathrm{Irr}(D_G^\ell(\tilde A)\mid  \psi_{A+x I} ).
	\end{equation}
	 For each $x\in \cO_{\ell'},$ $A+x I$ is cyclic. Let $\tilde x \in \cO_{2\ell'+1}$ be a lift of $x\in \cO_{\ell'}.$ Note that $\tilde A +\tilde x  I$  is a lift of  $A+x I,$ 
	 $D_G^\ell(\tilde A +\tilde x  I)= D_G^\ell(\tilde A)$  and $C_G^{\ell'}(\tilde {A}+ \tilde{x} I)=C_G^{\ell'}(\tilde A).$ Therefore  by Lemma~\ref{lem:odd r- construction-GL}, we obtain that each representation in $\mathrm{Irr}(D_G^\ell(\tilde A)\mid  \psi_{A+x I} )$ is one-dimensional and stabilised by $C_G^{\ell'}(\tilde A).$ Hence (3) follows by (\ref{eqn:abcds}).
\end{proof}

\begin{lemma}\label{lem:g^2 is idendity}
	For every $g\in D_S^{\ell'}(\tilde A),$ the matrix $g^2$ is in $ D_S^{\ell}(\tilde A).$
\end{lemma}
\begin{proof}
	Let $g\in  D_S^{\ell'}(\tilde A).$ Then by definition of $D_S^{\ell'}(\tilde A),$ there exists  $X \in C_{\GL_2(\cO_{2\ell'+1})} (\tilde{A})\cap M^1$ and $B \in M^{\ell'},$ such that $g=XB.$ Therefore $g^2=XBXB=X^2(X^{-1}BXB) \in D_S^{\ell}(\tilde A).$ The last inclusion follows because $X^2 \in  C_{\GL_2(\cO_{2\ell'+1})} (\tilde{A})\cap M^1$ and  $X^{-1}BXB \in M^\ell$ along with the fact that $g^2$ has determinant one.
\end{proof}
Note that by Lemma~\ref{lem:g^2 is idendity}, the quotient group 
$D_S^{\ell'}(\tilde A)/D_S^\ell(\tilde A) $  is  elementary abelian-2-group. 
For $\phi \in  \mathrm{Irr}(D_S^\ell(\tilde A)\mid  \psi_{[A]}),$ consider the alternating bilinear form $h'_\phi : D_S^{\ell'}(\tilde A)/D_S^{\ell}(\tilde A)\times D_S^{\ell'}(\tilde A)/D_S^{\ell}(\tilde A) \rightarrow 
\mathbb{C}^\times $ defined by $h'_\phi (xD_S^{\ell}(\tilde A) , y D_S^{\ell}(\tilde A))=\phi([x,y]).$ Note that $h'_\phi$ is well defined by Lemma~\ref{lem:D_S L phi}(2).
The following lemma gives a sufficient condition for the 
non-degeneracy of  the alternating bilinear form  $h'_\phi.$ 
\begin{lemma}\label{lem:h'-phi is non deg.}
	For $\phi \in  \mathrm{Irr}(D_S^\ell(\tilde A)\mid  \psi_{[A]}),$ assume the alternating bilinear form  $h'_\phi$ as above. If $\det(D_G^{\ell}(\tilde A))=\det(D_G^{\ell'}(\tilde A)),$ then $h'_\phi$ is non-degenerate.
\end{lemma}
\begin{proof}
	Let $\tilde \phi$ be an extension of $\phi$ into  $D_G^\ell(\tilde A)$ (obtained by Lemma~\ref{lem:D_S L phi}(3)). Then $\tilde \phi \in \mathrm{Irr}(D_G^{\ell}(\tilde A)\mid \psi_{A+v I})$ for some $v\in \cO_{\ell'}.$ 
	%\hfoot{notation "x" is also used for other things here.} 
	Consider the new  alternating bilinear form 
	$h_{\tilde \phi} : D_G^{\ell'}(\tilde A)/D_G^{\ell}(\tilde A)\times D_G^{\ell'}(\tilde A)/D_G^{\ell}(\tilde A)  \rightarrow \mathbb{C}^\times $ defined by $h_{\tilde \phi} (xD_G^{\ell}(\tilde A) , y D_G^{\ell}(\tilde A))=\tilde \phi([x,y]).$ Then by  Lemma~\ref{lem:non deg. of h phi}, we have $h_{\tilde \phi} $ is non-degenerate.
	
	To show $h'_\phi$ is non-degenerate, let $xD_S^{\ell}(\tilde A)\in D_S^{\ell'}(\tilde A)/D_S^{\ell}(\tilde A)$ be such that
	$xD_S^{\ell}(\tilde A)\neq D_S^{\ell}(\tilde A).$ Then clearly $x \in D_S^{\ell'}(\tilde A)\setminus D_S^{\ell}(\tilde A) \subseteq D_G^{\ell'}(\tilde A)\setminus D_G^{\ell}(\tilde A)$ and hence $ xD_G^{\ell}(\tilde A)\neq D_G^{\ell}(\tilde A).$ Therefore, since $h_{\tilde \phi} $ is non-degenerate,  there exists $y\in D_G^{\ell'}(\tilde A)$ such that $\tilde\phi([x,y])\neq 1.$ 
	By hypothesis (i.e., $\det(D_G^{\ell}(\tilde A))=\det(D_G^{\ell'}(\tilde A))$),  there exists $z\in D_G^{\ell}(\tilde A)$ such that $zy \in D_S^{\ell'}(\tilde A).$  Therefore  we have  the following. 
	$$\phi([x,zy])=\tilde\phi([x,zy])= \tilde\phi(z^{-1}[x,zy]z)=\tilde\phi([z^{-1},x])\tilde\phi([x,y])=\tilde\phi([x,y])\neq 1.$$
	Here the third equality follows from the fact that  $z^{-1}[x,zy]z=[z^{-1},x][x,y].$ Similarly the last equality follows because $C_G^{\ell'}(\tilde A)$ stabilizes  $\tilde \phi$ (by Lemma~\ref{lem:D_S L phi}(3)). So we have $h'_\phi(xD_G^{\ell}(\tilde A),zyD_G^{\ell}(\tilde A))=\phi([x,zy])\neq 1.$ Hence $h'_\phi$ is non-degenerate.
\end{proof}
%
%
%To prove Corollary~\ref{cor:beta inv construction char =0}, we need the following lemma, which is also used to prove Theorem~\ref{thm:construction-SL-r odd- char=2}.
\begin{lemma}\label{lem:image_of_det_of_DSL}
	Let  $\cO \in \dvrtwo$ and  
	$A \in M_2(\cO_{\ell'})$ be cyclic.
	\begin{enumerate}
		\item If $\mathrm{trace}(A)\in \cO_{\ell'}^\times,$ then $\det(D_G^{\ell}(\tilde A))=\det(D_G^{\ell'}(\tilde A)).$
		\item If $\cO \in \dvrtwozero$ and $\ell'\geq 2\mathrm{e}+1,$ then $\det(D_G^{\ell}(\tilde A))=\det(D_G^{\ell'}(\tilde A)).$
	\end{enumerate} 
\end{lemma}
\begin{proof}
		Note that $ D_G^{\ell'}(\tilde A)= D_G^{\ell}(\tilde A) M^{\ell'}$ and hence $\det(D_G^{\ell'}(\tilde A))= \det(D_G^{\ell}(\tilde A))  (1+\pi^{\ell'} \cO_{2\ell'+1}).$ Therefore to show  $\det(D_G^{\ell}(\tilde A))=\det(D_G^{\ell'}(\tilde A)),$ it is enough to show $1+\pi^{\ell'} \cO_{2\ell'+1} \subseteq \det(D_G^{\ell}(\tilde A)). $
	
	 To show (1), assume $\mathrm{trace}(A)\in \cO_{\ell'}^\times.$ For $y\in \cO_{2\ell'+1},$ observe that $ I +\pi^{\ell'} y \tilde{A} \in D_G^{\ell}(\tilde A)$ and 
	 $\det ( I +\pi^{\ell'} y \tilde{A})=1+\pi^{\ell'} ( \mathrm{trace}(\tilde A)y+  \det(\tilde A) \pi^{\ell'} y^2 ) .$
	 Let $z\in \cO_{2\ell'+1}.$ 
	 %Since $\mathrm{trace}(\tilde A)\in \cO_{2\ell'+1}^\times,$ 
	 Note that  $y= \mathrm{trace}(\tilde A)^{-1}z$ is a solution for $\mathrm{trace}(\tilde A)y+  \det(\tilde A) \pi^{\ell'} y^2=z$ in the residue field.
	 Since $\mathrm{trace}(\tilde A)\in \cO_{2\ell'+1}^\times,$ by Hensel's lemma, $\mathrm{trace}(\tilde A)y+  \det(\tilde A) \pi^{\ell'} y^2=z$ has a solution in $\cO_{2\ell'+1}.$
	 Therefore $1+\pi^{\ell'} \cO_{2\ell'+1}= \det(\{  I +\pi^{\ell'} y \tilde{A} \mid y \in \cO_{2\ell'+1} \})  \subseteq \det(D_G^{\ell}(\tilde A)). $
	 % Hence $\det(\{  I +\pi y \tilde{A} \mid y \in \cO_{2\ell'+1} \})=1+\pi \cO_{2\ell'+1}.$

	 To show (2), assume $\cO \in \dvrtwozero$ and $\ell'\geq 2\mathrm{e}+1.$
	  Let $w \in  \cO_{2\ell'+1}^\times$ be such that $2=\pi^{\mathrm{e}}w.$ For $y\in \cO_{2\ell'+1},$ observe that $(1+\pi^{\ell'-\mathrm{e}}y) I\in D_G^{\ell}(\tilde A) $ and $$\det((1+\pi^{\ell'-\mathrm{e}}y) I)=(1+\pi^{\ell'-\mathrm{e}}y)^2=1+2\pi^{\ell'-\mathrm{e}}y+\pi^{2\ell'-2\mathrm{e}}y^2=1+\pi^{\ell'}(wy+\pi^{\ell'-2\mathrm{e}}y^{2}).$$
	  Let $z \in \cO_{2\ell'+1}.$ 
	  Since $\ell'-2\mathrm{e}\geq 1,$  $y=w^{-1}z$ is a solution for $ wy+\pi^{\ell'-2\mathrm{e}}y^2=z$ in the residue field. Since $w \in  \cO_{2\ell'+1}^\times,$ by Hensel's lemma, $ wy+\pi^{\ell'-2\mathrm{e}}y^2=z$ has a solution in  $\cO_{2\ell'+1}.$ 
	   Therefore $1+\pi^{\ell'} \cO_{2\ell'+1}= \det(\{(1+\pi^{\ell'-\mathrm{e}}y) I \mid y \in \cO_{2\ell'+1} \})  \subseteq \det(D_G^{\ell}(\tilde A)). $

\end{proof}

%The following lemma will be useful in the following sections.
\begin{lemma}\label{lem:eqn-card}
	Let  $A \in M_2(\cO_{\ell'})$ be cyclic and $\tilde{A}   \in M_2(\cO_{2\ell'+1})$  be a lift  of $A.$ Then
	\begin{equation*}%\label{eqn-card}
	\frac{|D_S^{\ell'}(\tilde A)|}{|D_S^\ell(\tilde A)|}
	%=\frac{|D_G^{\ell'}(\tilde A)|}{|D_G^\ell(\tilde A)|}\times \frac{|\det(D_G^{\ell}(\tilde A))|}{|\det(D_G^{\ell'}(\tilde A))|} 
	=q^2 \times \frac{|\det(D_G^{\ell}(\tilde A))|}{|\det(D_G^{\ell'}(\tilde A))|}. 
	\end{equation*}
	
\end{lemma}
To prove Lemma~\ref{lem:eqn-card}, we need the following well known result. We also use this  in Section~\ref{sec:proof of construction cases 3}.
For it's proof, see for example \cite[Corollary~3.7]{MR1334228}. 
\begin{lemma} %\cite[Lemma 3.1(2)]{mypaper1}
	\label{lem:centralizer-form}
	For any $m \in \mathbb N$ and a cyclic matrix $A \in M_2(\cO_m),$ the  centralizer of $A$ in $\GL_2(\cO_m),$ denoted $C_{\GL_2(\cO_m)}(A) ,$ consists of invertible matrices of the form $x I + yA$ for $x, y \in \cO_m.$
\end{lemma}

\begin{proof}[Proof of Lemma~\ref{lem:eqn-card}]
	Note that $|D_S^{\ell'}(\tilde A)|=|D_G^{\ell'}(\tilde A)|/|\det(D_G^{\ell'}(\tilde A))|$ and $|D_S^{\ell}(\tilde A)|=|D_G^{\ell}(\tilde A)|/|\det(D_G^{\ell}(\tilde A))|.$ Therefore 
	$$\frac{|D_S^{\ell'}(\tilde A)|}{|D_S^\ell(\tilde A)|}
	=\frac{|D_G^{\ell'}(\tilde A)|}{|D_G^\ell(\tilde A)|}\times \frac{|\det(D_G^{\ell}(\tilde A))|}{|\det(D_G^{\ell'}(\tilde A))|} .$$
	Hence to show our result, it is enough to show $|D_G^{\ell'}(\tilde A)|/|D_G^\ell(\tilde A)|=q^2.$ But this can be obtained by direct calculations (using Lemma~\ref{lem:centralizer-form}).

\end{proof}
\begin{remark}\label{rmk:dl'/dl}
	Let  $A \in M_2(\cO_{\ell'})$ be cyclic and $\tilde{A}   \in M_2(\cO_{2\ell'+1})$  be a lift  of $A.$ Then $|D_G^{\ell'}(\tilde A)|/|D_G^\ell(\tilde A)|=q^2.$
\end{remark}
\begin{proposition}\label{prop:odd r- construction-GL*} %\pooja{Lemma 2.2.13(3), thesis}
 %For every cyclic $A \in M_2(\cO_{\ell'}),$ let $\tilde{A} \in M_2(\cO_{2\ell'+1})$ be a lift of $A.$ 
 	Let  $A \in M_2(\cO_{\ell'})$ be cyclic and $\tilde{A}   \in M_2(\cO_{2\ell'+1})$  be a lift  of $A.$
	Then  every $\rho \in \mathrm{Irr}(C_G^{\ell'}(\tilde A)\mid  \psi_{A}) $ has dimension $q.$
\end{proposition}
To prove Proposition~\ref{prop:odd r- construction-GL*}, we need the following result.
\begin{lemma}\label{lem:ext from DGL to CGL} 
	Let  $A \in M_2(\cO_{\ell'})$ be cyclic and $\tilde{A}   \in M_2(\cO_{2\ell'+1})$  be a lift  of $A.$
	Then  every $\theta \in \mathrm{Irr}(D_G^{\ell'}(\tilde A)\mid  \psi_{A}) $ has an extension to $C_G^{\ell'}(\tilde A).$
\end{lemma}
\begin{proof}
	 This lemma directly follows from a general result proved by Stasinski-Stevens,
	 % see \cite[Section~4]{MR3743488},
	  see\cite[Theorem~4.10]{MR3743488}.
\end{proof}
\begin{proof}[Proof of Proposition~\ref{prop:odd r- construction-GL*}]
	We first prove that  every  $\theta \in  \mathrm{Irr}(D_G^{\ell'}(\tilde A)\mid  \psi_{A}) $ has dimension $q.$ Let $\theta \in  \mathrm{Irr}(D_G^{\ell'}(\tilde A)\mid  \psi_{A}) .$ Choose $\chi \in \mathrm{Irr}(D_G^{\ell}(\tilde A)\mid  \psi_{A})$ such that $\theta \in  \mathrm{Irr}(D_G^{\ell'}(\tilde A)\mid  \chi) .$ Therefore, Lemma~\ref{lem:non deg. of h phi}, Proposition~\ref{prop:bilinear non deg.} and  Remark~\ref{rmk:dl'/dl} together give that  $\dim(\theta)=q.$ 
	
Since  $C_{\GL_2(\cO_{2\ell'+1})} (\tilde{A})$ is abelian, the quotient group $C_G^{\ell'}(\tilde A) / D_G^{\ell'}(\tilde A)$ must be abelian. Therefore by Clifford theory,  to show the proposition, it is enough to show that each  $\theta \in  \mathrm{Irr}(D_G^{\ell'}(\tilde A)\mid  \psi_{A}) $ extends to $C_G^{\ell'}(\tilde A).$ But the later follows from Lemma~\ref{lem:ext from DGL to CGL}.

\end{proof}
The following corollary of Proposition~\ref{prop:odd r- construction-GL*} will be used in this section and 
Section~\ref{Sec:group-algebras}.
\begin{corollary}\label{cor:dim q [A]}
	Let  $A \in M_2(\cO_{\ell'})$ be cyclic and $\tilde{A}   \in M_2(\cO_{2\ell'+1})$  be a lift  of $A.$
	Every $\chi \in \mathrm{Irr}(C_G^{\ell'}(\tilde A)\mid  \psi_{[A]} )$ has dimension $q.$
\end{corollary}
\begin{proof}
	Let  $\chi \in \mathrm{Irr}(C_G^{\ell'}(\tilde A)\mid  \psi_{[A]} ).$
	Note that $\mathrm{Irr}(M^\ell\mid  \psi_{[A]} )=\{\psi_{A+x I} \mid x\in \cO_{\ell'} \}.$ Therefore
	\begin{equation}\label{eqn:cor}
		\mathrm{Irr}(C_G^{\ell'}(\tilde A)\mid  \psi_{[A]} )=\bigcup_{x\in \cO_{\ell'}}\mathrm{Irr}(C_G^{\ell'}(\tilde A)\mid  \psi_{A+x I} ).
	\end{equation}
	For each $x\in \cO_{\ell'},$ $A+x I$ is cyclic. Let $\tilde x \in \cO_{2\ell'+1}$ be a lift of $x\in \cO_{\ell'}.$ Then we have $\tilde A +\tilde x  I$  is a lift of  $A+x I$ and   $C_G^{\ell'}(\tilde {A}+ \tilde{x} I)=C_G^{\ell'}(\tilde A).$
	Therefore by Proposition~\ref{prop:odd r- construction-GL*}, we  obtain that for each $x\in \cO_{\ell'}, $ every representation in $\mathrm{Irr}(C_G^{\ell'}(\tilde A)\mid  \psi_{A+x I} )$ has dimension $q.$ Hence from  (\ref{eqn:cor}), we have $\dim(\chi)=q.$
\end{proof}

 %%%%%%%%%%%%%%%%%%%%%%%%%%%%%%%
Next we prove Theorem~\ref{thm:construction-SL-r odd- char=0 and 2(half)}.
\begin{proof}[Proof of Theorem~\ref{thm:construction-SL-r odd- char=0 and 2(half)}]
	Note that (1) of Theorem~\ref{thm:construction-SL-r odd- char=0 and 2(half)} follows from Lemma~\ref{lem:D_S L phi}(1). 
By Lemma~\ref{lem:image_of_det_of_DSL}, we have 
$\det(D_G^{\ell}(\tilde A))=\det(D_G^{\ell'}(\tilde A)).$ Therefore (2) of Theorem~\ref{thm:construction-SL-r odd- char=0 and 2(half)} follows from Proposition~\ref{prop:bilinear non deg.} and Lemmas~\ref{lem:h'-phi is non deg.} and \ref{lem:eqn-card}.   
%Now we prove (3) of Theorem~\ref{thm:construction-SL-r odd- char=0 and 2(half)}. 
By Corollary~\ref{cor:dim q [A]}, we  obtain that each  $\chi \in \mathrm{Irr}(C_G^{\ell'}(\tilde A)\mid  \psi_{[A]} )$ has dimension $q.$ Thus (3) of Theorem~\ref{thm:construction-SL-r odd- char=0 and 2(half)} follows because $C_S^{\ell'}(\tilde A)\leq C_G^{\ell'}(\tilde A)$ and $\hat \phi$ has dimension $q.$
	
	Next we show (4) of Theorem~\ref{thm:construction-SL-r odd- char=0 and 2(half)}. We first claim that the representation $\mathrm{Ind}_{C_S^{\ell'}(\tilde A)}^{ \SL_2(\cO_{2\ell'+1})}(\rho)$ is irreducible for every $\rho \in \mathrm{Irr}(C_S^{\ell'}(\tilde A)\mid \psi_{[A]}).$ Then by Clifford
	theory, we have 
	$$\mathrm{Irr}(\SL_2(\cO_{2 \ell'+1})\mid  \psi_{[A]})= \{\mathrm{Ind}_{C_S^{\ell'}(\tilde A)}^{ \SL_2(\cO_{2\ell'+1})}(\rho) \mid \rho \in \mathrm{Irr}(C_S^{\ell'}(\tilde A)\mid \psi_{[A]}) \} .$$
	 Hence to show (4) of Theorem~\ref{thm:construction-SL-r odd- char=0 and 2(half)}, it is enough to prove that any $\rho_1 , \rho_2 \in \mathrm{Irr}(C_S^{\ell'}(\tilde A)\mid \psi_{[A]})$ satisfy $\mathrm{Ind}_{C_S^{\ell'}(\tilde A)}^{ \SL_2(\cO_{2\ell'+1})}(\rho_1) \cong \mathrm{Ind}_{C_S^{\ell'}(\tilde A)}^{ \SL_2(\cO_{2\ell'+1})}(\rho_2)$ if and only if $\rho_1 \cong \rho_2^g$ for some $g\in  C_{\SL_2(\cO_{2\ell'+1})} (\psi_{[A]}). $
	 This follows because of the Clifford theory and the fact that $C_S^{\ell'}(\tilde A)$ is
	 a normal subgroup of $ C_{\SL_2(\cO_{2\ell'+1})} (\psi_{[A]}). $

Now we show the claim. Let $\rho \in \mathrm{Irr}(C_S^{\ell'}(\tilde A)\mid \psi_{[A]}).$ By Clifford theory, it is enough to show that  the stabilizer of $\rho$ in  $ C_{\SL_2(\cO_{2\ell'+1})} (\psi_{[A]})$  is $C_S^{\ell'}(\tilde A).$
By (1)-(3) of  Theorem~\ref{thm:construction-SL-r odd- char=0 and 2(half)} and Clifford theory, there exists one-dimensional representation $\phi \in \mathrm{Irr}(D_S^\ell(\tilde A)\mid \psi_{[A]}) $ such that $\rho$ is an extension of the unique representation $\hat{\phi} \in  \mathrm{Irr}(D_S^{\ell'}(\tilde A)\mid  \phi).$ 
%
 %let  $\rho \in \mathrm{Irr}(C_S^{\ell'}(\tilde A)\mid \hat{\phi}).$ We claim that  the stabilizer of $\rho$ in  $ C_{\SL_2(\cO_{2\ell'+1})} (\psi_{[A]})$  is $C_S^{\ell'}(\tilde A).$ Then clearly (3) follows  from %Theorem~\ref{clifford theory}(2)	Clifford theory.To show the claim, 
 We recall from Lemma~\ref{lem:stabilizer-form} that 
	$$ C_{\SL_2(\cO_{2\ell'+1})} (\psi_{[A]})= C_S^{\ell'}(\tilde A)\{e_\lambda=\mat 1{\tilde{a}^{-1}\lambda}01  \mid \lambda \in \mathrm{h}_{\tilde{A}}^{\ell'} \}.$$
	Note that $e_\lambda \in  C_S^{\ell'}(\tilde A)$ for all $\lambda \in \pi^{\ell'}\cO_{2\ell'+1}.$ 
	Consider the set $L_\rho:=\{ \lambda \in  \mathrm{h}_{\tilde{A}}^{\ell'}\mid e_\lambda \,\, \mathrm{stabilizes} \, \, \rho  \}.$ 
	To show 
	the stabilizer of $\rho$ in  $ C_{\SL_2(\cO_{2\ell'+1})} (\psi_{[A]})$  is $C_S^{\ell'}(\tilde A),$
 it is enough  to show that $L_\rho\subseteq\pi^{\ell'}\cO_{2\ell'+1}.$ Since $\rho $ is an extension of $\hat\phi$ and $\hat\phi$ is the unique irreducible representation lying above $\phi,$ we must have 
	\begin{equation}\label{eqn:incl 1}
		L_\rho \subseteq \{ \lambda \in  \mathrm{h}_{\tilde{A}}^{\ell'}\mid e_\lambda \,\, \mathrm{stabilizes} \, \, \phi \}.
	\end{equation}
	
	For $\lambda \in \mathrm{h}_{\tilde{A}}^{\ell'},$ if $e_\lambda$ stabilizes $\phi,$ then
	the group $D_S^\ell(\tilde A) \langle e_\lambda  \rangle  $ stabilises $\phi$ and  the quotient group $D_S^\ell(\tilde A) \langle e_\lambda  \rangle / D_S^\ell(\tilde A)$ is cyclic.  Thus by \cite[Corollary~11.22]{MR2270898},  
	$\phi$ extends to $D_S^\ell(\tilde A) \langle e_\lambda  \rangle  .$ 
%	
%	
%	For $\lambda \in \mathrm{h}_{\tilde{A}}^{\ell'},$ if $e_\lambda$ stabilizes $\phi,$ then it is easy to see that $\phi([x,y])=1$ for all $x,y\in D_S^{\ell}(\tilde A)\langle e_\lambda  \rangle .$
%	Therefore by Lemma~\ref{prop:extn-commutater-relation}, $\phi$ extends to $D_S^{\ell}(\tilde A)\langle e_\lambda  \rangle .$
%	%
	Hence by definition of $E^\prime_{\tilde A},$ we obtain that 
	\begin{equation}\label{eqn:incl 2}
		\{ \lambda \in  \mathrm{h}_{\tilde{A}}^{\ell'}\mid e_\lambda \,\, \mathrm{stabilizes} \, \, \phi \} \subseteq E^\prime_{\tilde A}. 
	\end{equation}
	
	For the case $\cO \in \dvrtwozero$
	and $\ell' \geq 2\mathrm{e}+1,$ note that  $r=2\ell'+1\geq 4\mathrm{e}+3$ and hence  by Theorem~\ref{S_A-in-number},
	we have  $E^\prime_{\tilde A}=\pi^{\ell'}\cO_{2\ell'+1}.$ Similarly, for the case $\cO \in \dvrtwoplus,$
	$\ell' \geq 1$   and  $\mathrm{trace}(A)\in \cO_{\ell'}^\times,$ by Proposition~\ref{prop:char 2 E*A}, we have 
	$E^\prime_{\tilde A}=\pi^{\ell'}\cO_{2\ell'+1}.$ 
	Therefore by combining (\ref{eqn:incl 1}) and (\ref{eqn:incl 2}) we obtain $L_\rho\subseteq\pi^{\ell'}\cO_{2\ell'+1}.$ Hence the claim holds.
	
%	Next assume  $\cO \in \dvrtwoplus,$
%	 $\ell' \geq 1$   and  $\mathrm{trace}(A)\in \cO_{\ell'}^\times$. Then by Proposition~\ref{prop:char 2 E*A}(2), we have 
%	  $E^\prime_{\tilde A}=\pi^{\ell'}\cO_{2\ell'+1}.$ Therefore by combining (\ref{eqn:incl 1}) and (\ref{eqn:incl 2}) we obtain $L_\rho\subseteq\pi^{\ell'}\cO_{2\ell'+1}.$ Hence the claim holds.
%	
%	
%	Since $\ell' \geq 2\mathrm{e}+1,$  $r=2\ell'+1\geq 4\mathrm{e}+3$ and hence  by Theorem~\ref{S_A-in-number},
%	we have  $E^\prime_{\tilde A}=\pi^{\ell'}\cO_{2\ell'+1}.$ Therefore by combining (\ref{eqn:incl 1}) and (\ref{eqn:incl 2}) we obtain $L_\rho\subseteq\pi^{\ell'}\cO_{2\ell'+1}.$ Hence the claim holds.

\end{proof}	
\begin{corollary}\label{cor:beta inv construction char =0}
	Let  $\cO \in \dvrtwozero$ and $A \in M_2(\cO_{\ell'})$ be cyclic such that $\mathrm{trace}(A)\in \cO_{\ell'}^\times.$ Then Theorem~\ref{thm:construction-SL-r odd- char=0 and 2(half)} holds for $\ell'>\mathrm{e}.$
\end{corollary}
%To prove Corollary~\ref{cor:beta inv construction char =0}, we need the following lemma, which is also used to prove Theorem~\ref{thm:construction-SL-r odd- char=2}.
%\begin{lemma}%\label{lem:image_of_det_of_DSL}
%	Let  $\cO \in \dvrtwo$ and  $A \in M_2(\cO_{\ell'})$ be cyclic such that $\mathrm{trace}(A)\in \cO_{\ell'}^\times.$ Then $\det(D_G^{\ell}(\tilde A))=\det(D_G^{\ell'}(\tilde A)).$
%\end{lemma}
%\begin{proof}
%	Note that $\{  I +\pi y \tilde{A} \mid y \in \cO_{2\ell'+1} \}\subseteq D_G^{\ell}(\tilde A)\subseteq D_G^{\ell'}(\tilde A)$ and $\det(D_G^{\ell'}(\tilde A)) \subseteq 1+\pi \cO_{2\ell'+1}.$ Therefore to show  $\det(D_G^{\ell}(\tilde A))=\det(D_G^{\ell'}(\tilde A)),$ it is enough to show  $\det(\{  I +\pi y \tilde{A} \mid y \in \cO_{2\ell'+1} \})=1+\pi \cO_{2\ell'+1}.$ Note that $\det ( I +\pi y \tilde{A})=1+\pi ( \mathrm{trace}(\tilde A)y+  \det(\tilde A) \pi y^2 ) .$ Let $z\in \cO_{2\ell'+1}.$ Since $\mathrm{trace}(\tilde A)\in \cO_{2\ell'+1}^\times,$ observe that  $y= \mathrm{trace}(\tilde A)^{-1}z$ is a solution for $\mathrm{trace}(\tilde A)y+  \det(\tilde A) \pi y^2=z$ in the residue field.
%	Therefore by Hensel's lemma, $\mathrm{trace}(\tilde A)y+  \det(\tilde A) \pi y^2=z$ has a solution in $\cO_{2\ell'+1}.$ Hence $\det(\{  I +\pi y \tilde{A} \mid y \in \cO_{2\ell'+1} \})=1+\pi \cO_{2\ell'+1}.$
%\end{proof}
\begin{proof}
	Let $\mathrm{trace}(A)\in \cO_{\ell'}^\times$ and $\ell'>\mathrm{e}.$
	Note that in the proof of Theorem~\ref{thm:construction-SL-r odd- char=0 and 2(half)},  we used the assumption $\ell'\geq 2\mathrm{e}+1$ only to prove  $\det(D_G^{\ell}(\tilde A))=\det(D_G^{\ell'}(\tilde A))$ and  $E^\prime_{\tilde A}=\pi^{\ell'}\cO_{2\ell'+1}.$ 
	So if we show $\det(D_G^{\ell}(\tilde A))=\det(D_G^{\ell'}(\tilde A))$ and  $E^\prime_{\tilde A}=\pi^{\ell'}\cO_{2\ell'+1},$ then Theorem~\ref{thm:construction-SL-r odd- char=0 and 2(half)} holds in this case. Note that Lemma~\ref{lem:image_of_det_of_DSL}(1)  and Theorem~\ref{S_A-in-number}(1)  give  $\det(D_G^{\ell}(\tilde A))=\det(D_G^{\ell'}(\tilde A))$ and $E^\prime_{\tilde A}=\pi^{\ell'}\cO_{2\ell'+1}$ respectively.
\end{proof}

%\begin{remark}\label{rmk:thm 0---2}
%	Note that the arguments used to prove Theorem~\ref{thm:construction-SL-r odd- char=0} are independent of $\Char(\cO),$ except those used to show $\det(D_G^{\ell}(\tilde A))=\det(D_G^{\ell'}(\tilde A))$ and  $E^\prime_{\tilde A}=\pi^{\ell'}\cO_{2\ell'+1}.$
%\end{remark}
%
%
%
%
%%
%Next we prove Theorem~\ref{thm:construction-SL-r odd- char=2}.
%\begin{proof}[Proof of Theorem~\ref{thm:construction-SL-r odd- char=2}]
%	Since $\mathrm{trace}(A)\in \cO_{\ell'}^\times,$ by Lemma~\ref{lem:image_of_det_of_DSL} and Proposition~\ref{prop:char 2 E*A}(2), we have 
% $\det(D_G^{\ell}(\tilde A))=\det(D_G^{\ell'}(\tilde A))$ and  $E^\prime_{\tilde A}=\pi^{\ell'}\cO_{2\ell'+1}.$ Now the theorem follows by the same arguments used to prove Theorem~\ref{thm:construction-SL-r odd- char=0} (see Remark~\ref{rmk:thm 0---2}). 
%\end{proof}
%

\section{Proof of Theorem~\ref{thm:construction-SL-r odd- char=2 (non inv trace)}}\label{sec:proof of construction cases 3}

In this section we prove Theorem~\ref{thm:construction-SL-r odd- char=2 (non inv trace)}.
Through out this section we assume  $\cO \in \dvrtwoplus,$ 
%and $r\geq 2$ is odd. $r=2\ell'+1.$  
 we fix  $A = \mat 0{a^{-1} \alpha}a{\beta} \in M_2(\cO_{\ell'})$ and its lift $\tilde{A} = \mat 0{\tilde{a}^{-1}\tilde \alpha}{\tilde{a}}{\tilde \beta} \in M_2(\cO_{2\ell'+1}).$ We also assume $\mathrm{trace}(A)=\beta \in  \pi\cO_{\ell'}.$ 
%
%Recall the definition  of $\bol s$ and $\bol k$.  \hfoot{ckeck defn of k} 
%\begin{eqnarray}
%\bol  k & = & \min \{ \val({\tilde{\beta}}), \ell' \}. \nonumber \\
%\bol s & = &\begin{cases}
%2 \lfloor \bol k/2 \rfloor + 1 &,\,\, \mathrm{if}\,\,   \alpha = v^2 \,\,\mathrm{mod}\,(\pi^ \bol k )  \\ m &,\,\, \mathrm{if}\,\, \alpha = v_1^2 + \pi^m v_2^2 \,\,\mathrm{mod}\,(\pi^ \bol k ) \,\,\mathrm{for}\,\,\mathrm{odd}\,\, m < \bol k \,\,\mathrm{and}\,\, v_2 \in \cO_{\ell'}^\times. 
%\end{cases}\nonumber
%\end{eqnarray}
%Also for $\lambda \in \cO_{2\ell'+1},$ recall the following notations.
%\begin{eqnarray}
%\bol i_\lambda & = &\val({\lambda}). \nonumber \\
%\bol  j_\lambda & = & \min \{ \val({\lambda+ \tilde{\beta}}), \ell' \}. \nonumber 
%\end{eqnarray}
%
The following lemmas will be useful in this section.
\begin{lemma}\label{lem:E'tilde A = E tilde A}
		%Let  $\cO \in \dvrtwoplus$ and $\ell' \geq 1.$  For a cyclic $A= \mat 0 {a^{-1}\alpha } a\beta \in M_2(\cO_{\ell'})$ such that $\mathrm{trace}(A)=\beta\in \pi \cO_{\ell'},$ we have 
Assume $A\in M_2(\cO_{\ell'})$ as above. Then 	$E^\prime_{\tilde{A}}=E_{\tilde{A}}.$
\end{lemma}
\begin{proof}
	By Lemma~\ref{lem:link lemma}, we have $E^\prime_{\tilde{A}}\setminus\pi^{\ell'}\cO_{2 \ell'+1} \subseteq \mathrm{h}_{\tilde{A}}^{\ell}.$ Since $\tilde{\beta}\in  \pi \cO_{2\ell'+1},$ for every $\lambda \in \pi^{\ell'}\cO_{2 \ell'+1},$ we have $\lambda (\lambda +\tilde{ \beta} ) = 0  \,\,\mathrm{mod}\,(\pi^{\ell} ).$ Therefore, by definition of $ \mathrm{h}_{\tilde{A}}^{\ell},$ we obtain $\pi^{\ell'}\cO_{2 \ell'+1} \subseteq  \mathrm{h}_{\tilde{A}}^{\ell}.$ Hence we have $E^\prime_{\tilde{A}}\subseteq \mathrm{h}_{\tilde{A}}^{\ell}.$ Now the lemma follows because $E_{\tilde{A}}=E^\prime_{\tilde{A}}\cap  \mathrm{h}_{\tilde{A}}^{\ell}$ (by Corollary~\ref{cor:EtildeA and E'tildeA}).
\end{proof}
\begin{remark}\label{rmk:E_tilde A in pi Or}
	Since $\tilde \beta\in \pi\cO_{2 \ell'+1},$ for each $\lambda\in \cO_{2 \ell'+1}^\times ,$ $\lambda(\lambda+\tilde \beta)\in  \cO_{2 \ell'+1}^\times .$ Therefore by definition of $\mathrm{h}_{\tilde{A}}^{\ell'},$ we must have $\mathrm{h}_{\tilde{A}}^{\ell'}\subseteq \pi \cO_{2 \ell'+1}.$ In particular, $E^\prime_{\tilde{A}}\subseteq \pi \cO_{2 \ell'+1}.$
\end{remark}
%\begin{lemma}\label{lem:quotient and stab}
%	\hfoot{why $C_S^{\ell'}(\tilde A)\mathbb E^\prime_{\tilde{A}}$ is a group?}
% The quotient group $C_S^{\ell'}(\tilde A)\mathbb E^\prime_{\tilde{A}}/ D_S^{\ell}(\tilde A)$ is abelian.
%		
%\end{lemma}

\begin{lemma}\label{lem:quotient and stab}
	For each $g \in C_S^{\ell'}(\tilde A)\mathbb E^\prime_{\tilde{A}}, $ 
the element $g^2 \in D_S^{\ell}(\tilde A).$ 
	
\end{lemma}

 \begin{proof}%\hfoot{simplify this proof.}
 %	Note that for a group $G,$ if $z^2=1$ for all $z\in G,$ then $G$ is abelian.
	Let $X \in C_{\GL_2(\cO_{2\ell'+1})} (\tilde{A})$ and $B= I+\pi^{\ell'}B' \in M^{\ell'}$ be such that $XB \in C_S^{\ell'}(\tilde A).$ Let $e_\lambda\in \mathbb E^\prime_{\tilde{A}}.$
	To show the result, 
	it is enough to 
	prove that $(XBe_\lambda)^2 \in D_S^{\ell}(\tilde A).$
	 Since  $\det (XBe_\lambda)=1,$ it is enough to show that  $(XBe_\lambda)^2 \in D_G^{\ell}(\tilde A).$ i.e., we have to show $(XBe_\lambda)^2 \tilde A=\tilde A (XBe_\lambda)^2 \,\, \mathrm{mod}\,(\pi^{\ell'+1}).$ 
	 	Note that 
	 \begin{eqnarray*}
	 	(XBe_\lambda)^2&=&X( I+\pi^{\ell'}B')e_\lambda X( I+\pi^{\ell'}B')e_\lambda\\
	 	&=& Xe_\lambda Xe_\lambda +\pi^{\ell'}(X B'e_\lambda X e_\lambda +X e_\lambda X B' e_\lambda)+\pi^{2\ell'}X B'e_\lambda X B'e_\lambda .
	 \end{eqnarray*}
	 Since 
	 $e_\lambda= I \,\, \mathrm{mod}\,(\pi)$ (by Remark~\ref{rmk:E_tilde A in pi Or}, $\lambda \in \pi\cO_{2 \ell'+1}$),  we must have 
	 \begin{equation}\label{eqn:l1}
	 (XBe_\lambda)^2= Xe_\lambda Xe_\lambda +\pi^{\ell'}(X B' X  +X  X B' )\,\, \mathrm{mod}\,(\pi^{\ell'+1}).
	 \end{equation}
	
	By Lemma~\ref{lem:centralizer-form}, there exists $x,y \in \cO_{2\ell'+1}$ such that $X=x I + y \tilde{A}.$
	 Since   $ \lambda \in E^\prime_{\tilde{A}}=E_{\tilde{A}} \subseteq \mathrm{h}_{\tilde{A}}^{\ell}$ (by Lemma~\ref{lem:E'tilde A = E tilde A}), we have $\lambda(\lambda + \tilde{ \beta})=0 \,\, \mathrm{mod}\, (\pi^{\ell'+1} ).$ Also note that $e_\lambda = e_\lambda^{-1}.$ Therefore by \cite[Lemma~5.1(1)]{mypaper1} we have $e_\lambda \tilde{A} e_\lambda= \tilde{A} + \lambda I \,\, \mathrm{mod}\,(\pi^{\ell'+1}).$ Hence we obtain that 
	$Xe_\lambda Xe_\lambda = (x I + y \tilde{A})(x I + y \tilde{A}+ \lambda y I)\,\, \mathrm{mod}\, (\pi^{\ell'+1} ),$ which 
 implies 
 	 \begin{equation}\label{eqn:l2}
  (Xe_\lambda Xe_\lambda) \tilde{A}=\tilde{A} (Xe_\lambda Xe_\lambda) \,\, \mathrm{mod}\, (\pi^{\ell'+1} ).
\end{equation}

	%The last equality follows because $\tilde{A}^2= \alpha  I + \tilde{ \beta} \tilde{A}.$
	
	Now we observe that 
	$$ \pi^{\ell'}(X B' X  +X  X B' )=\pi^{\ell'}X( B' X  +  X B' )=\pi^{\ell'}y(x I + y \tilde{A}) (B'  \tilde{A}  +  \tilde{A} B' ).$$
		Since $\tilde A=\mat{0}{\tilde{a}^{-1}\tilde \alpha }{\tilde{a}}{0} \,\, \mathrm{mod}\,(\pi),$
		% and $X=X^{-1}=\mat{x}{\tilde{a}^{-1}\tilde \alpha y}{\tilde{a}y}{x} \,\, \mathrm{mod}\,(\pi),$ 
		 for $B'=[b_{ij}],$ we have 
		$$ B'  \tilde{A}  +  \tilde{A} B' =\mat{\tilde{a} b_{12}+ \tilde{a}^{-1} \tilde \alpha b_{21}}{ \tilde{a}^{-1} \tilde \alpha( b_{11}+  b_{22})}{ \tilde{a}( b_{11}+ b_{22})}{\tilde{a} b_{12}+ \tilde{a}^{-1} \tilde \alpha b_{21}}= (\tilde{a} b_{12}+ \tilde{a}^{-1} \tilde \alpha b_{21}) I + ( b_{11}+ b_{22}) \tilde{A}   \,\, \mathrm{mod}\,(\pi).$$
		 Therefore 
		 	$ \pi^{\ell'}(X B' X  +X  X B' )=\pi^{\ell'}y(x I + y \tilde{A}) ((\tilde{a} b_{12}+ \tilde{a}^{-1} \tilde \alpha b_{21}) I + ( b_{11}+ b_{22}) \tilde{A} ) \,\, \mathrm{mod}\,(\pi^{\ell'+1}),$ which implies 
		 	 \begin{equation*}%\label{eqn:l3}
		 		 \pi^{\ell'}(X B' X  +X  X B' )\tilde A=\tilde A  \pi^{\ell'}(X B' X  +X  X B' ) \,\, \mathrm{mod}\,(\pi^{\ell'+1}).
		 	\end{equation*}
		 This along with (\ref{eqn:l1}) and (\ref{eqn:l2})  give 
		 $(XBe_\lambda)^2 \tilde A=\tilde A (XBe_\lambda)^2 \,\, \mathrm{mod}\,(\pi^{\ell'+1}).$ 	Hence $(XBe_\lambda)^2 \in D_G^{\ell}(\tilde A).$

\end{proof}

\begin{proposition}\label{prop:abelian qoutient}
	Let $G$ be a finite group and  $N$ be a normal subgroup of $G$ such that $G/N$ is abelian. Suppose $\chi: N\rightarrow \mathbb{C}^\times$ is a one-dimensional representation such that it is stabilised by $G.$ Then 	there exists a group $ M_\chi$ such that 
	\begin{enumerate}
		\item $N \leq  M_\chi \leq G. $ 
		\item $\chi$ extends to $  M_\chi.$  Further each $\rho \in \mathrm{Irr}(M_\chi\mid \chi)$ is an extension of $\chi.$
		\item  For each $\rho \in \mathrm{Irr}( M_\chi\mid \chi),$ the induced representation $\mathrm{Ind}_{M_\chi}^{G}(\rho)$ is irreducible.

	\end{enumerate}
	
\end{proposition}

To prove Proposition~\ref{prop:abelian qoutient}, we need the following lemmas.
\begin{lemma} (Diamond Lemma)%\cite[Lemma~5.4]{MR2684153}
	\label{diamond-lemma} 
	Let $H$ and $K$ be two subgroups of a group $G$ such that $H$ is normal in $G$ and $G = HK.$ Let $\chi_1$ and $\chi_2$ be one-dimensional representations of $H$ and $K$ respectively such that $\chi_1(khk^{-1}) = \chi_1(h)$ for all $h \in H$ and $k \in K,$ and $\chi_1|_{H \cap K } = \chi_2|_{H\cap K}.$ Then there exists a unique one-dimensional representation of $G$ extending both $\chi_1$ and $\chi_2$ simultaneously.  
\end{lemma}
\begin{proof} For proof, see \cite[Lemma~5.4]{MR2684153}. 
\end{proof}

\begin{lemma}\label{lem:extn-commutater-relation}
	Let $N$ be a normal subgroup of a finite group $G.$ Suppose $\chi$ is a one-dimensional representation of $N.$ Then $\chi$ has an extension to $G$ if and only if $[G,G] \cap N$ is contained in the kernel of $\chi.$
\end{lemma}

\begin{proof}
	If $\chi$ has an extension to $G,$ then clearly  $[G,G] \cap N$ is contained in the kernel of $\chi.$ To prove the backward implication, assume $[G,G] \cap N$ is contained in the kernel of $\chi.$ 
	Note that for $g\in G$ and $n\in N,$ the element $gng^{-1}n^{-1}\in[G,G] \cap N $ and therefore $\chi(gng^{-1}n^{-1})=1.$ So we obtain that  $G$ stabilises $\chi,$ in particular
	the subgroup $[G,G]$ stabilises $\chi.$ Therefore by Lemma~\ref{diamond-lemma} and the hypothesis, $\chi$ extends to a one-dimensional representation $\hat{\chi}$ of $N[G,G]$ such that restriction of  $\hat{\chi}$ to $[G,G]$  is trivial. Thus we can consider $\hat{\chi}$ as  a representation of the subgroup $\{n[G,G] \mid n \in N\}$ of $G/[G,G].$ Since $G/[G,G]$ is abelian, $\hat{\chi}$ extends to $G/[G,G].$ Let $\tilde\chi$ be such an extension. 
	Then the one-dimensional representation $\tilde\chi 
	\comp 
	\theta : G \rightarrow 
	\mathbb{C}^\times,$ 
	where $\theta:G \rightarrow G/[G,G]$ is the canonical projection,  is an extension of $\chi.$
\end{proof}

\begin{proof}[Proof of Proposition~\ref{prop:abelian qoutient}]
	Consider the set of subgroups $\mathcal{F}=\{M : N\leq M\leq G  \text{ and } [M,M] \subseteq \mathrm{ker}(\chi)\}.$ Note that $N\in \mathcal{F}$ and hence $\mathcal{F}\neq \emptyset.$ Since $G$ is finite, the set $\mathcal{F}$ must be finite.   Therefore there exists a maximal subgroup $M_\chi \in \mathcal{F}.$ i.e., for $M\in \mathcal{F}$ if $M_\chi \leq M,$ then $M=M_\chi.$ We show that this maximal subgroup $M_\chi$ satisfies the required conditions.
	
	 Clearly $N \leq  M_\chi \leq G. $ Since $[M_\chi,M_\chi] \subseteq \mathrm{ker}(\chi),$ by Lemma~\ref{lem:extn-commutater-relation}, $\chi$ extends to $  M_\chi.$ Since  the quotient group $M_\chi/N$ is abelian (because $G/N$ is abelian),  by Clifford theory,  each $\rho \in \mathrm{Irr}(M_\chi\mid \chi)$ is an extension of $\chi.$ Let $\rho \in \mathrm{Irr}(M_\chi\mid \chi).$ To show the induced representation $\mathrm{Ind}_{M_\chi}^{G}(\rho)$ is irreducible, by Clifford theory, it is enough to show that the stabiliser  $C_G(\rho)=M_\chi.$ Note that $M_\chi \subseteq C_G(\rho).$ 
	 Suppose on the contrary that
	 $C_G(\rho)\neq M_\chi.$ Then there exists $x\in C_G(\rho)\setminus M_\chi.$
	 Consider the subgroup $M_{\chi}^\prime$   of $G$ generated by $M_\chi$ and $x.$ We claim that $M_{\chi}^\prime \in \mathcal{F}.$ Assume the claim for now. Since $M_\chi$ is a maximal subgroup in $\mathcal{F}$ and $M_{\chi}\subseteq M_{\chi}^\prime,$ we must have $M_{\chi}^\prime=M_{\chi}.$ It is a contradiction to the fact that $x\notin M_\chi.$ Hence we must have $C_G(\rho)= M_\chi.$    
	 
	 To show the claim, note that since $M_\chi \trianglelefteq G,$ we have $M_{\chi}^\prime =\{zx^i : z\in M_\chi \text{ and }i \in \Z\} .$
	 For $zx^i, z'x^j \in M_{\chi}^\prime,$ we have $[zx^i, z'x^j]=zx^i z'x^jx^{-i}z^{-1}x^{-j} z^{\prime-1}=z(x^i z'x^{-i})(x^jz^{-1}x^{-j}) z^{\prime-1}.$ Therefore 
	 \begin{eqnarray*}
	\chi([zx^i, z'x^j])&=&\rho([zx^i, z'x^j])\\
	&=&\rho(z)\rho(x^i z'x^{-i})\rho(x^jz^{-1}x^{-j}) \rho(z^{\prime-1})\\
	&=&\rho(z)\rho( z')\rho(z^{-1}) \rho(z^{\prime-1})\\
	&=&1.
	 \end{eqnarray*}
 Here the third equality follows because $x\in C_G(\rho).$ Hence we must have $[M_{\chi}^\prime,M_{\chi}^\prime] \subseteq \mathrm{ker}(\chi).$ 
	This together with  $N\leq M_{\chi}^\prime \leq G$ implies $M_{\chi}^\prime \in \mathcal{F}.$ Hence the claim holds.
	 
\end{proof}

Next we prove Theorem~\ref{thm:construction-SL-r odd- char=2 (non inv trace)}.
\begin{proof}[Proof of Theorem~\ref{thm:construction-SL-r odd- char=2 (non inv trace)}]
	Note that (1) of Theorem~\ref{thm:construction-SL-r odd- char=2 (non inv trace)} follows from Lemma~\ref{lem:D_S L phi}(1). To show (2) of Theorem~\ref{thm:construction-SL-r odd- char=2 (non inv trace)}, let $\phi \in \mathrm{Irr}(D_S^\ell(\tilde A)\mid \psi_{[A]}). $ We first show that the stabiliser $\mathrm{Stab}(\phi)$ of $\phi$ in  $C_{\SL_2(\cO_{2\ell'+1})} (\psi_{[A]})$ is contained in $C_S^{\ell'}(\tilde A)\mathbb E^\prime_{\tilde{A}}.$
	Since $C_{\SL_2(\cO_{2\ell'+1})} (\psi_{[A]})=C_S^{\ell'}(\tilde A)\mathrm{H}_{\tilde{A}}^{\ell'}$ (by Lemma~\ref{lem:stabilizer-form}), it is enough to show that 
	$ \mathrm{H}_{\tilde{A}}^{\ell'} \cap \mathrm{Stab}(\phi) \subseteq \mathbb E^\prime_{\tilde{A}}.$ For each $e_\lambda\in \mathrm{H}_{\tilde{A}}^{\ell'} \cap \mathrm{Stab}(\phi),$
	the group $D_S^\ell(\tilde A) \langle e_\lambda  \rangle  $ stabilises $\phi$ and  the quotient group $D_S^\ell(\tilde A) \langle e_\lambda  \rangle / D_S^\ell(\tilde A)$ is cyclic.  Thus by \cite[Corollary~11.22]{MR2270898},  
	  $\phi$ extends to $D_S^\ell(\tilde A) \langle e_\lambda  \rangle  .$ Therefore by the definition of 
	  $\mathbb{E}^\prime_{\tilde{A}},$  we obtain
	   $e_\lambda \in  \mathbb{E}^\prime_{\tilde{A}}. $  Hence 
	  $\mathrm{H}_{\tilde{A}}^{\ell'} \cap \mathrm{Stab}(\phi)\subseteq \mathbb{E}^\prime_{\tilde{A}}. $

	 Since  $\mathrm{Stab}(\phi)\subseteq C_S^{\ell'}(\tilde A)\mathbb E^\prime_{\tilde{A}},$ by Lemma~\ref{lem:quotient and stab},
	 $g^2D_S^\ell(\tilde A)=D_S^\ell(\tilde A)$ for all $g\in \mathrm{Stab}(\phi). $ Therefore 
	  the quotient group $\mathrm{Stab}(\phi)/D_S^\ell(\tilde A)$ must be abelian. Hence by Proposition~\ref{prop:abelian qoutient},
	there exists $D_S^\ell(\tilde A) \leq  \mathbb M_\phi \leq \mathrm{Stab}(\phi)$  such that 
	\begin{enumerate}
		\item $\phi$ extends to $ \mathbb M_\phi$  and each $\rho \in \mathrm{Irr}(\mathbb M_\phi\mid \phi)$ is an extension of $\phi.$
		\item  For each $\rho \in \mathrm{Irr}(\mathbb M_\phi\mid \phi),$ the induced representation $\mathrm{Ind}_{\mathbb M_\phi}^{\mathrm{Stab}(\phi)}(\rho)$ is irreducible.

	\end{enumerate}
	   Clearly this group $\mathbb M_\phi$ satisfies  (2)(a) and (2)(b) of Theorem~\ref{thm:construction-SL-r odd- char=2 (non inv trace)}. To show  $\mathbb M_\phi$ satisfies  (2)(c) of Theorem~\ref{thm:construction-SL-r odd- char=2 (non inv trace)}, let $\rho \in \mathrm{Irr}(\mathbb M_\phi\mid \phi).$ 
	   Note that by (2), $\mathrm{Ind}_{\mathbb M_\phi}^{\mathrm{Stab}(\phi)}(\rho)\in \mathrm{Irr}(\mathrm{Stab}(\phi)\mid \phi).$ Therefore,  by Clifford theory, the induced representation 
	   $$\mathrm{Ind}_{\mathbb M_\phi}^{C_{\SL_2(\cO_{2\ell'+1})} (\psi_{[A]})}(\rho)\cong
	    \mathrm{Ind}_{\mathrm{Stab}(\phi)}^{C_{\SL_2(\cO_{2\ell'+1})} (\psi_{[A]})}(\mathrm{Ind}_{\mathbb M_\phi}^{\mathrm{Stab}(\phi)} (\rho))$$
	     is irreducible. Hence, again by Clifford theory, we obtain that the induced representation 
	     $\mathrm{Ind}_{\mathbb M_\phi}^{\SL_2(\cO_{2\ell'+1})}(\rho)$
	     is irreducible. So  $\mathbb M_\phi$ satisfies  (2)(c) of Theorem~\ref{thm:construction-SL-r odd- char=2 (non inv trace)}. 

\end{proof}

%Next we prove  that if $2\bol k - \bol s \geq \ell,$ then the group $C_S^{\ell}(\tilde A)\mathbb E^\prime_{\tilde{A}} $ work  as $\mathbb M_\phi$ for all $\phi \in \mathrm{Irr}(D_S^\ell(\tilde A)\mid \psi_{[A]}). $ In particular, we have the following result.
%\begin{corollary}
%		Let  $\cO \in \dvrtwoplus$
%	and $\ell' \geq 1.$  The following hold for every $A=\mat 0 {a^{-1}\alpha} a \beta \in M_2(\cO_{\ell'})$ such that $2\bol k - \bol s \geq \ell.$
%	% and  it's lift $\tilde{A} =  \mat 0 {\tilde{a}^{-1} \tilde{\alpha}}{\tilde a}{\tilde{\beta}}  \in M_2(\cO_{r})$ of  $A.$  such that $\mathrm{trace}(A)\in \pi \cO_{\ell'}.$
%
%	\begin{enumerate}
%		\item Every representation $\rho \in \mathrm{Irr}(C_S^{\ell}(\tilde A)\mathbb E^\prime_{\tilde{A}} \mid \psi_{[A]}) $  is one-dimensional.
%		\item  For each $\rho \in \mathrm{Irr}(C_S^{\ell}(\tilde A)\mathbb E^\prime_{\tilde{A}} \mid \psi_{[A]}),$ the induced representation $\mathrm{Ind}_{C_S^{\ell}(\tilde A)\mathbb E^\prime_{\tilde{A}}}^{ \SL_2(\cO_{2\ell'+1})}(\rho)$ is irreducible.
%	\end{enumerate}
%	
%	
%\end{corollary}
%
%	
%

\section{Group Algebra of $\SL_2(\cO_{r})$}
	%{Proof of Theorem~\ref{thm:group-algebras}}
	\label{Sec:group-algebras}
In this section, we prove  Theorem~\ref{thm:group-algebras}. 
	For a finite group $G,$ the polynomial 
	\[
	\calp_G(X) := \sum_{\rho \in \mathrm{Irr}(G) } X^{\dim(\rho)}
	\]
	is called the {\it representation zeta polynomial} of $G.$ As a consequence result of Wedderburn's theorem (see \cite[Theorem~14, Chapter~15]{MR1138725}), 
	%it is easy to see that   
	for finite groups $G_1$ and $G_2,$ the group algebras $\mathbb{C}[G_1]$ and $ \mathbb{C}[G_2]$ are isomorphic if and only if $\calp_{G_1}(X) = \calp_{G_2}(X).$
Therefore to prove  Theorem~\ref{thm:group-algebras}, we prove the following result. 
\begin{thm}\label{thm:zeta poly not eaqal- gp algebra} Let $\cO \in \dvrtwozero$ with ramification index $\ee $ and $\cO' \in \dvrtwoplus$ such that $\cO/\wp \cong \cO'/\wp'.$  Then  
	\[
	\calp_{\SL_2(\cO_{r})}(X)\neq \calp_{\SL_2(\cO'_{r}) }(X)
	\]
	for any $r \geq 2 \mathrm{e}+2.$
\end{thm}

We recall the following lemmas from \cite{mypaper1}.
%When $\bar{A}$ is cyclic, the following lemma is easy and well known in the literature.
\begin{lemma}\cite[Lemma~3.3]{mypaper1}
	\label{gl-centralizer-cardinality}
	Let $A \in M_2(\cO_m)$ such that $\bar{A}$ is cyclic. Then  
	\[
	|C_{\GL_2(\cO_{m})}(A)| = \begin{cases} (q^2-1)q^{2 m -2} &
	\,\, \mathrm{for} \,\, \bar{A} \,\, \mathrm{irreducible}, \\ (q-1)^2 q^{2m -2}&
	\,\, \mathrm{for} \,\, \bar{A} \,\, \mathrm{split\,\, semisimple}, \\ (q^2-q)q^{2 m -2} & \,\, \mathrm{for} \,\, 
	\bar{A} \,\, \mathrm{ split, non-semisimple} . \end{cases}
	\]
\end{lemma}
\begin{lemma} \cite[Lemma~3.4]{mypaper1}
	\label{lem:image-det-map}
	For 
	$A = \mat 0 {a^{-1}\alpha} a \beta \in M_2(\cO_m)$ such that $\beta$ is invertible, the group $C_{\GL_2(\cO_{m})}(A)$ maps onto $\cO_{m}^\times$ under the determinant map. 
	
\end{lemma} 

\begin{lemma}\label{lem:card-[A]-stab*}
	Let $A \in M_2(\cO_{\ell'})$ be  cyclic.
	% such that  $\bar{A}$ split  non-semisimple. 
	Then 
	$$|C_{\SL_2(\cO_{r})} (\psi_{[A]})|=\begin{cases}
		(q-1)2^{m_1}& ,\,\,\mathrm{if} \,\, \bar{A} \,\, \mathrm{split\,\, semisimple} \\ 	(q+1)2^{m_2}&,\,\, \mathrm{if} \,\, \bar{A} \,\, \mathrm{irreducible} \\ 2^{m_3}&,\,\, \mathrm{if} \,\, \bar{A} \,\, \mathrm{split\,\, non-semisimple},
	\end{cases}$$
	for some $m_1, m_2, m_3 \geq 0.$
\end{lemma}
 \begin{proof}
 	From  Lemma~\ref{lem:stabilizer-form}, we have $C_{\SL_2(\cO_{r})} (\psi_{[A]})=C_{\SL_2(\cO_{r})} (\psi_{A})\mathrm{H}_{\tilde{A}}^{\ell'}=C_S^{\ell'}(\tilde A)\mathrm{H}_{\tilde{A}}^{\ell'}.$ %where $\mathrm{H}_{\tilde{A}}^{\ell'} = \{ \mat 1{\tilde{a}^{-1}x}01 \mid x \in \cO_{r } \, \, \mathrm{with}\,\,  2x = 0 \,\, \mathrm{mod}\,   (\pi^{\ell'} ) ,\,\, x(x+\tilde \beta) = 0 \,\, \mathrm{mod}\,   (\pi^{\ell'} )  \}.$
 	By definition of $\mathrm{H}_{\tilde{A}}^{\ell'},$ it is easy to observe that $\mathrm{H}_{\tilde{A}}^{\ell'}$ is a $2\text{-group}.$ Therefore to show the result, it is enough to show  that
 	$$|C_S^{\ell'}(\tilde A)|=\begin{cases}
 		(q-1)2^{m'_1}& ,\,\,\mathrm{if} \,\, \bar{A} \,\, \mathrm{split\,\, semisimple} \\ 	(q+1)2^{m'_2}&,\,\, \mathrm{if} \,\, \bar{A} \,\, \mathrm{irreducible} \\ 2^{m'_3}&,\,\, \mathrm{if} \,\, \bar{A} \,\, \mathrm{split\,\, non-semisimple},
 	\end{cases}$$
 	for some $m'_1, m'_2, m'_3 \geq 0.$
 	But by definition of $C_S^{\ell'}(\tilde{A}),$   we obtain the following.
 	\begin{eqnarray*}
 		|C_S^{\ell'}(\tilde{A})|&=&	|C_{\SL_2(\cO_{\ell'})}(A)| \times |K^{\ell'}|\\
 		&=& \frac{ |C_{\GL_2(\cO_{\ell'})}(A)|}{|\det(C_{\GL_2(\cO_{\ell'})}(A))|}\times q^{3\ell} \\
 		&=&
 		\begin{cases}
 			(q-1)q^{3 \ell+\ell' -1}& ,\,\,\mathrm{if} \,\, \bar{A} \,\, \mathrm{split\,\, semisimple} \\ 	(q+1)q^{3 \ell+\ell' -1}&,\,\, \mathrm{if} \,\, \bar{A} \,\, \mathrm{irreducible} \\   \frac{(q-1)q^{3 \ell +2\ell'-1}}{|\det(C_{\GL_2(\cO_{\ell'})}(A))|}&,\,\, \mathrm{if} \,\, \bar{A} \,\, \mathrm{split\,\, non-semisimple}.
 		\end{cases}
 		%\times q^{3(\ell'+1)} 
 	\end{eqnarray*}
 	The last equality follows from Lemmas~\ref{gl-centralizer-cardinality} and \ref{lem:image-det-map} along with the fact that $\mathrm{trace}(A)\in \cO_{\ell'}^\times$ if $\bar{A}$ is either  split  semisimple or irreducible. Since $q$ is a power of $2,$ we are done for the cases  $\bar{A}$ split  semisimple and  $\bar{A}$  irreducible.
 	
 	For the case $\bar{A}$ split non-semisimple,  to show $|C_S^{\ell'}(\tilde A)|=2^{m'_3}$ for some $m'_3\geq 1,$ we need to show that  $|\det(C_{\GL_2(\cO_{\ell'})}(A))|=(q-1)2^t$ for some $t\geq0.$ The later is equivalent to show   $(q-1)$ divides $|\det(C_{\GL_2(\cO_{\ell'})}(A))|,$  because $|\det(C_{\GL_2(\cO_{\ell'})}(A))|$ divides $|\cO_{\ell'}^\times|=(q-1)q^{\ell'-1}.$  
 	Note that $\{x^2=\det(x I) \mid x \in \cO_{\ell'}^\times\}\subseteq \det(C_{\GL_2(\cO_{\ell'})}(A)).$ Therefore the natural projection $\det(C_{\GL_2(\cO_{\ell'})}(A))\rightarrow \mathbb{F}_q^\times$ is onto and hence $(q-1)=|\mathbb{F}_q^\times|$ divides $|\det(C_{\GL_2(\cO_{\ell'})}(A))|.$ This completes  proof of the lemma.
 \end{proof}

 \begin{lemma}\label{lem:dim divisibility}
 	Let $A \in M_2(\cO_{\ell'})$ be a cyclic matrix. Let 
 	$$\gamma_A=\begin{cases}
 		q+1& ,\,\,\mathrm{if} \,\, \bar{A} \,\, \mathrm{split\,\, semisimple} \\ 	q-1&,\,\, \mathrm{if} \,\, \bar{A} \,\, \mathrm{irreducible} \\ 	q^2-1&,\,\, \mathrm{if} \,\, \bar{A} \,\, \mathrm{split\,\, non-semisimple}. \end{cases}$$
 	Then for any 
 	$\rho \in \mathrm{Irr}( \SL_2(\cO_{r}) \mid \psi_{[A]}),$  the term $\gamma_A$ divides $\dim(\rho).$
 \end{lemma}
 \begin{proof}
 	Let $A\in  M_2(\cO_{\ell'})$ be  cyclic and $\rho \in \mathrm{Irr}( \SL_2(\cO_{r}) \mid \psi_{[A]}).$ Then by Clifford theory, there exists $\theta\in \mathrm{Irr}(C_{\SL_2(\cO_{r})} (\psi_{[A]})\mid \psi_{[A]})$ such that $\rho\cong  \mathrm{Ind}_{C_{\SL_2(\cO_{r})} (\psi_{[A]})}^{ \SL_2(\cO_{r})}(\theta).$ Therefore $\frac{|\SL_2(\cO_{r})|}{|C_{\SL_2(\cO_{r})} (\psi_{[A]})|}$ divides $\dim(\rho).$ Hence to show $\gamma_A$  divides $\dim(\rho),$ it is enough to show that $\gamma_A$  divides $\frac{|\SL_2(\cO_{r})|}{|C_{\SL_2(\cO_{r})} (\psi_{[A]})|}.$ But the later follows from Lemma~\ref{lem:card-[A]-stab*} along with  $|\SL_2(\cO_{r})|=(q^2 -1)q^{3r-2}.$
 	%$(q^2-1)$ divides $|\SL_2(\cO_{r})|.$ 

 \end{proof}
 
 \begin{lemma}\label{lem:irred. split non sem  dim} 
 	%\hfoot{Is this statement fine?}
 	Let $A \in M_2(\cO_{\ell'})$ be cyclic such that $\bar{A}$ is irreducible (resp. split  non-semisimple).
 	Then for any 
 	$\rho \in \mathrm{Irr}( \SL_2(\cO_{r}) \mid \psi_{[A]}),$ $\dim(\rho)\leq(q-1)q^{r-1}$ (resp. $\dim(\rho)\leq(q^2-1)q^{r-2}$).
 \end{lemma}
 \begin{proof}
 	Let $A \in M_2(\cO_{\ell'})$ be cyclic such that $\bar{A}$ is either irreducible or split  non-semisimple. Define 
 	$$d_A=\begin{cases}
 		(q-1)q^{r-1}& \mathrm{if} \,\, \bar{A} \, \, \mathrm{is} \,\, \mathrm{irreducible}\\
 		(q^2 -1)q^{r-2}& \mathrm{if} \,\, \bar{A} \, \, \mathrm{is} \,\, \mathrm{split}\,\, \mathrm{non-semisimple}.
 	\end{cases}$$
 	Let $\rho \in \mathrm{Irr}( \SL_2(\cO_{r}) \mid \psi_{[A]}).$
 	Note that by Clifford theory, there exists $\chi\in \mathrm{Irr}(C_S^{\ell'}(\tilde A)\mid \psi_{[A]})$ such that $\rho$ is an irreducible constituent of the induced representation $ \mathrm{Ind}_{C_S^{\ell'}(\tilde A)}^{ \SL_2(\cO_{r})}(\chi).$ Therefore  to show $\dim(\rho)\leq
 	d_A,$ it is enough to show that $\dim(\mathrm{Ind}_{C_S^{\ell'}(\tilde A)}^{ \SL_2(\cO_{r})}(\chi))\leq d_A,$ which is equivalent to show that $ \frac{| \SL_2(\cO_{r})|}{|C_S^{\ell'}(\tilde A)|}\times \dim(\chi)\leq d_A.$ Note that $|\SL_2(\cO_{r})|=(q^2 -1)q^{3r-2}$ and 
 	$$|C_S^{\ell'}(\tilde{A})|=	|C_{\SL_2(\cO_{\ell'})}(A)| \times |K^{\ell'}|=\frac{ |C_{\GL_2(\cO_{\ell'})}(A)|}{|\det(C_{\GL_2(\cO_{\ell'})}(A))|}\times q^{3\ell} \geq \frac{ |C_{\GL_2(\cO_{\ell'})}(A)|}{q-1}\times q^{3\ell-\ell'+1}.$$
 	The last inequality follows because
 	% $|C_{\GL_2(\cO_{\ell'})}(A)|=(q^2 -1)q^{2\ell'-2}$ (form Lemma~\ref{gl-centralizer-cardinality}) and
 	$|\det(C_{\GL_2(\cO_{\ell'})}(A))|\leq |\cO_{\ell'}^\times|=(q-1)q^{\ell'-1}.$ %(since $\mathrm{trace}(A)\in \cO_{\ell'}^\times,$ from Lemma~\ref{lem:image-det-map}).
 	Therefore  we obtain  
 	\begin{eqnarray*}
 		\frac{| \SL_2(\cO_{r})|}{|C_S^{\ell'}(\tilde A)|}&\leq& \frac{(q^2 -1)q^{3r-2}\times (q-1)}{|C_{\GL_2(\cO_{\ell'})}(A)|\times q^{3\ell-\ell'+1}}\\
 		&=& \begin{cases}
 			(q-1)q^{2\ell'-1} & \mathrm{if} \,\, \bar{A} \, \, \mathrm{is} \,\, \mathrm{irreducible}\\
 			(q^2 -1)q^{2\ell'-2}& \mathrm{if} \,\, \bar{A} \, \, \mathrm{is} \,\, \mathrm{split}\,\, \mathrm{non-semisimple}.
 		\end{cases}.
 	\end{eqnarray*}
 	The last equality follows from Lemma~\ref{gl-centralizer-cardinality}.
 	Therefore to show $ \frac{| \SL_2(\cO_{r})|}{|C_S^{\ell'}(\tilde A)|}\times \dim(\chi)\leq d_A,$ it is enough to show that  $\dim(\chi)\leq q^{r-2\ell'}.$ For even $r,$ we have $\ell=\ell'$ and hence $ C_S^{\ell'}(\tilde A)=C_S^{\ell}(\tilde A).$ Therefore  by \cite[Lemma 3.2]{mypaper1}, 
 	%Lemma~\ref{lem:centralizer-sl}(3), 
 	 $\dim(\chi)=1=q^{r-2\ell'}.$ So we proved the result for even $r.$
 	
 	Next assume $r$ is odd. Note that $r-2\ell'=1$ in this case.
 	Since $\chi\in \mathrm{Irr}(C_S^{\ell'}(\tilde A)\mid \psi_{[A]})$ and $C_S^{\ell'}(\tilde A)\leq C_G^{\ell'}(\tilde A),$ we have $\dim(\chi)\leq \dim(\theta)$ for some $\theta \in\mathrm{Irr}(C_G^{\ell'}(\tilde A)\mid \psi_{[A]}).$ 
 From Corollary~\ref{cor:dim q [A]},  we have  $ \dim(\theta)=q$ for all $\theta \in\mathrm{Irr}(C_G^{\ell'}(\tilde A)\mid \psi_{[A]}).$ So we obtain $\dim(\chi)\leq q=q^{r-2\ell'}.$ Hence the result holds for odd $r$ also.
 	
 \end{proof}
 
To prove Theorem~\ref{thm:zeta poly not eaqal- gp algebra}, we need the following results.% will be very useful in this section.

\begin{lemma}\label{lem:odd r dims from cyclic}
		%Assume either $\ell'\geq1$ (if $\cO \in \dvrtwoplus$) or $\ell'>\mathrm{e}$ (if $\cO \in \dvrtwozero$).
	Assume either  $\cO \in \dvrtwoplus$ and $\ell'\geq1,$  or  $\cO \in \dvrtwozero$ and $\ell'>\mathrm{e}$.
	Let $A \in M_2(\cO_{\ell'})$ be a cyclic matrix. Then any $\rho \in \mathrm{Irr}( \SL_2(\cO_{2 \ell'+1}) \mid \psi_{[A]})$ satisfies the following.
	\begin{enumerate}
		\item $\dim(\rho)=\begin{cases}
		(q+1)q^{2\ell'}& ,\,\,\mathrm{if} \,\, \bar{A} \,\, \mathrm{split\,\, semisimple} \\ 	(q-1)q^{2\ell'}&,\,\, \mathrm{if} \,\, \bar{A} \,\, \mathrm{irreducible} . \end{cases}$
		\item For $\bar{A}$ split  non-semisimple, $(q^2-1)$ divides $\dim(\rho).$ Further $\dim(\rho)<
		(q+1)q^{2\ell'}.$ 
		%and {\color{red}*********************}\hfoot{.....}
		
	\end{enumerate}
\end{lemma}

\begin{proof}
	Let $\rho \in \mathrm{Irr}( \SL_2(\cO_{2 \ell'+1}) \mid \psi_{[A]}).$ To show (1),
	note that if $\bar{A}$ is irreducible or  split semisimple, then $\mathrm{trace}(A)\in \cO_{\ell'}^\times.$ Therefore  by Corollary~\ref{cor:beta inv construction char =0}  and  Theorem~\ref{thm:construction-SL-r odd- char=0 and 2(half)},  there exists a $q$-dimensional representation $\chi\in \mathrm{Irr}(C_S^{\ell'}(\tilde A)\mid \psi_{[A]})$ such that $\rho \cong \mathrm{Ind}_{C_S^{\ell'}(\tilde A)}^{ \SL_2(\cO_{2\ell'+1})}(\chi).$ Hence $\dim(\rho)=q\times \frac{| \SL_2(\cO_{2\ell'+1})|}{|C_S^{\ell'}(\tilde A)|}.$  By definition of $C_S^{\ell'}(\tilde{A}),$ we have the following equality (this holds for any cyclic $A\in M_2(\cO_{\ell'})$).
	\begin{equation}\label{C_S L' cardinality}
	|C_S^{\ell'}(\tilde{A})|=	|C_{\SL_2(\cO_{\ell'})}(A)| \times |K^{\ell'}|  = \frac{ |C_{\GL_2(\cO_{\ell'})}(A)|}{|\det(C_{\GL_2(\cO_{\ell'})}(A))|}\times q^{3(\ell'+1)} .
	\end{equation}
	Since $\mathrm{trace}(A)\in \cO_{\ell'}^\times,$ by Lemma~\ref{lem:image-det-map}, we have $ |\det(C_{\GL_2(\cO_{\ell'})}(A))|=|\cO_{\ell'}^\times|=(q-1)q^{\ell'-1}.$	Therefore 
	$$\dim(\rho)=q\times\frac{|\SL_2(\cO_{2\ell'+1})|\times (q-1)q^{\ell'-1}}{ |C_{\GL_2(\cO_{\ell'})}(A)|\times q^{3(\ell'+1)} }=  \frac{ (q^2 -1)(q-1)q^{4\ell'-2}}{ |C_{\GL_2(\cO_{\ell'})}(A)|}.$$
	The last inequality follows from the fact that  $|\SL_2(\cO_{2\ell'+1})|=(q^2 -1)q^{6\ell'+1}.$
	Now (1) follows from Lemma~\ref{gl-centralizer-cardinality}.

	To show (2), let $A\in  M_2(\cO_{\ell'})$ be  cyclic such that  $\bar{A}$ split  non-semisimple. Since $\rho \in \mathrm{Irr}( \SL_2(\cO_{2 \ell'+1}) \mid \psi_{[A]}),$ by Clifford theory, there exists $\theta\in \mathrm{Irr}(C_{\SL_2(\cO_{2\ell'+1})} (\psi_{[A]})\mid \psi_{[A]})$ such that $\rho\cong  \mathrm{Ind}_{C_{\SL_2(\cO_{2\ell'+1})} (\psi_{[A]})}^{ \SL_2(\cO_{2\ell'+1})}(\theta).$ Therefore $\frac{|\SL_2(\cO_{2\ell'+1})|}{|C_{\SL_2(\cO_{2\ell'+1})} (\psi_{[A]})|}$ divides $\dim(\rho).$ Hence to show $(q^2 -1)$  divides $\dim(\rho),$ it is enough to show that $(q^2 -1)$  divides $\frac{|\SL_2(\cO_{2\ell'+1})|}{|C_{\SL_2(\cO_{2\ell'+1})} (\psi_{[A]})|}.$ But the later follows because 
	$(q^2-1)$ divides $|\SL_2(\cO_{2\ell'+1})|$ and $|C_{\SL_2(\cO_{2\ell'+1})} (\psi_{[A]})|=2^m$ for some $m\geq 1$ (by Lemma~\ref{lem:card-[A]-stab*}).

	Next we show  $\dim(\rho)<
	(q+1)q^{2\ell'}.$ Since $\bar{A}$ is  split  non-semisimple, by Lemma~\ref{lem:irred. split non sem  dim}, we have $\dim(\rho)\leq (q^2 -1) q^{2\ell'-1}.$ Now the result follows because $(q^2 -1) q^{2\ell'-1} <	(q+1)q^{2\ell'}.$ Hence (2) holds.
\end{proof}
\begin{lemma}\label{lem:max dim from cyclic}
	For $r\geq 2,$ any irreducible representation of $\SL_2(\cO_{r})$ has dimension  less than or  equal to $(q+1) q^{r-1}.$
\end{lemma}

\begin{proof} 
	From Onn~\cite[Theorem~1.4]{MR2456275}, any irreducible representation of $\GL_2(\cO_{r})$ has dimension less than or equal to $(q+1) q^{r-1}.$ The lemma follows because  $\SL_2(\cO_{r})$ is a subgroup of $\GL_2(\cO_{r}).$ 
\end{proof}

\begin{lemma}\label{lem:second max dim from cyclic}
	Let $\ell>1.$ Then there does not exist an irreducible representation $\rho$ of $\SL_2(\cO_{2\ell-1})$ such that $\dim(\rho)=(q-1)q^{2\ell-1}.$
\end{lemma}
\begin{proof}
	Let $\rho \in \mathrm{Irr}(\SL_2(\cO_{2\ell-1})).$ 	
	Suppose on the contrary that $\dim(\rho)=(q-1)q^{2\ell-1}.$ Note that $\rho$ satisfies either of the following.
	\begin{enumerate}
		\item $\rho\in \mathrm{Irr}(\SL_2(\cO_{2\ell-1})\mid \psi_{[A]})$ for some cyclic $A\in M_2(\cO_{\ell-1}).$
		\item $\rho\in \mathrm{Irr}(\SL_2(\cO_{2\ell-1})/K^{2\ell-2}).$
		%\cong \mathrm{Irr}(\SL_2(\cO_{2\ell-2})).$
	
	\end{enumerate}
Observe that $\dim(\rho)=(q-1)q^{2\ell-1}$ is neither divisible by  $(q+1)$ nor by $(q^2 -1).$ Therefore
 if  $\rho$ satisfies (1),  then by Lemma~\ref{lem:dim divisibility},  $\bar{A}$ irreducible. So by Lemma~\ref{lem:irred. split non sem  dim}, we must have  $\dim(\rho)\leq (q-1)q^{2\ell-2}.$ It is a contraction to our hypothesis $\dim(\rho)=(q-1)q^{2\ell-1}.$
 So we assume that $\rho$ satisfies (2). Note that $ \SL_2(\cO_{2\ell-1})/K^{2\ell-2}\cong \SL_2(\cO_{2\ell-2}).$  Therefore by Lemma~\ref{lem:max dim from cyclic}, we must have $\dim(\rho)\leq (q+1)q^{2\ell-3}.$ It is a contradiction to our hypothesis $\dim(\rho)=(q-1)q^{2\ell-1}.$ Therefore the lemma holds.
\end{proof}

\begin{proposition}\label{prop:second degrre for odd r}
	Let $\cO \in \dvrtwo$ and $\ell' \geq 1.$ Then the coefficients of $X^{(q-1)q^{2\ell'-1}}$ in $	\calp_{\SL_2(\cO_{2 \ell'+1})}(X)$ and in $	\calp_{\SL_2(\cO_{2 \ell'})}(X)$ are  equal. 
\end{proposition}
\begin{proof}
	Consider the set of cyclic irreducible representations of $\SL_2(\cO_{2 \ell'+1}),$ that is the one that lie above $\psi_{[A]}$ for a cyclic $A.$ Recall from Section~\ref{sec: basic-framework}, this is exactly the set of primitive irreducible representations, denoted by $\mathrm{Irr}^\pr(\SL_2(\cO_{2 \ell'+1})),$ of $\SL_2(\cO_{2 \ell'+1}).$ The primitive representation zeta polynomial corresponding to the cyclic representations of $\SL_2(\cO_{2 \ell'+1})$ is defined as below.  
	\[
	\calp^\pr_{\SL_2(\cO_{2 \ell'+1})}(X) = \sum_{\rho \in \mathrm{Irr}^\pr(\SL_2(\cO_{2 \ell'+1})) } X^{\dim(\rho)}.
	\]
Note that in case $A$ is not cyclic then $\psi_{[A]} $ is a trivial one-dimensional representation of $K^{2\ell'}.$ Therefore  $\SL_2(\cO_{2 \ell'+1})$ has the following expression for the representation zeta polynomial. 
	\[
	\calp_{\SL_2(\cO_{2 \ell'+1})}(X)   = \calp^\pr_{\SL_2(\cO_{2 \ell'+1})}(X) + \calp_{\SL_2(\cO_{2\ell'})} (X).
	\]
	So, to show the result, it is enough to show that the coefficient of $X^{(q-1)q^{2\ell'-1}}$ in $\calp^\pr_{\SL_2(\cO_{2\ell'+1})} (X)$ is zero. For that we need to show    $\dim(\rho)\neq (q-1)q^{2\ell'-1},$  for all $\rho \in \mathrm{Irr}^\pr(\SL_2(\cO_{2 \ell'+1})).$ This follows because for each  $\rho \in \mathrm{Irr}^\pr(\SL_2(\cO_{2 \ell'+1})),$
	by Lemma~\ref{lem:odd r dims from cyclic}, we have $\dim(\rho)$ is either in $\{(q+1)q^{2\ell'},(q-1)q^{2\ell'}\}$ or divisible by $(q^2 -1).$

\end{proof} 

\begin{lemma}\label{lem:second high dim iff irreducible}
	Assume $\cO \in \dvrtwozero .$ Let $\ell>\mathrm{e}$ and $A \in M_2(\cO_\ell)$ be a cyclic matrix. 
	%Then the following are true. 
	\begin{enumerate} 
		\item Any $\rho \in \mathrm{Irr}( \SL_2(\cO_{2 \ell}) \mid \psi_{[A]})$ satisfies $\dim(\rho) = (q-1) q^{2 \ell-1}$ if and only if $\bar{A}$ is irreducible.
		\item The group $\SL_2(\cO_{2 \ell})$ has $\frac{(q^2-1) q^{2 \ell -2 }}{2}$ many inequivalent primitive irreducible representations of dimension $(q-1) q^{2 \ell-1}.$ 
	\end{enumerate} 
\end{lemma} 
\begin{proof}
	For (1),  the forward implication is  obtained by Lemma~\ref{lem:dim divisibility}. For the converse, let $A \in M_2(\cO_\ell)$ be such that $\bar{A}$ is irreducible. Note that without loss of generality we can take $A= \mat0{a^{-1}\alpha}a\beta.$ 
	Since $\bar{A}$ is irreducible, $\beta=\mathrm{trace}(A)$ is invertible. Therefore by definition of $\mathrm{h}_{\tilde{A}}^\ell$ along with our assumption $\ell>\mathrm{e},$ we have 
	$\mathrm{h}_{\tilde{A}}^\ell =\pi^\ell\cO_{2\ell} .$ Hence $C_{\SL_2(\cO_{2\ell})} (\psi_{[A]})= C_{\SL_2(\cO_{2\ell})}(\psi_{A})$ (by Lemma~\ref{lem:stabilizer-form}(3)). 
	Therefore for $\rho \in \mathrm{Irr}( \SL_2(\cO_{2 \ell}) \mid \psi_{[A]}),$ by 
	Clifford theory there exists 
	$\phi \in   \mathrm{Irr}( C_{\SL_2(\cO_{2\ell})}(\psi_{A}) \mid \psi_{[A]})$ such that $\rho \cong \mathrm{Ind}_{C_{\SL_2(\cO_{2\ell})}(\psi_{A})}^{ \SL_2(\cO_{2\ell})}(\phi) .$ Since $C_{\SL_2(\cO_{2\ell})}(\psi_{A})=C_S^\ell(\tilde{A}),$ by \cite[Lemma~3.2]{mypaper1},
	$\phi$ is one-dimensional. Therefore 
	$\dim(\rho) = \frac{|\SL_2(\cO_{2\ell})|}{|C_{\SL_2(\cO_{2\ell})}(\psi_{A})|}= (q-1) q^{2 \ell-1}.$ The last equality follows from  Lemmas~\ref{lem:image-det-map} and  \ref{gl-centralizer-cardinality}.
	
	Next, for (2) we note that from  \cite[Theorem~2.1]{mypaper1}, %Lemma~\ref{lem:number-of-orbits-split semisimple},
	 the number of orbits of $\psi_{[A]}$ with $\bar{A}$ is irreducible is $ \frac{(q-1)(q^{\ell-1})}{2}.$ Therefore
	from (1), it is enough to show that for $A \in M_2(\cO_\ell)$ with $\bar{A}$ is irreducible, $| \mathrm{Irr}( \SL_2(\cO_{2 \ell}) \mid \psi_{[A]})|=(q+1)q^{\ell-1}.$ Let $A= \mat0{a^{-1}\alpha}a\beta \in M_2(\cO_\ell)$ be such that $\bar{A}$ is irreducible. Note that we have  $C_{\SL_2(\cO_{2\ell})} (\psi_{[A]})=C_{\SL_2(\cO_{2\ell})}(\psi_{A}).$ Now the result follows as below by Clifford theory %\hfoot{any ref. ?}.
	$$| \mathrm{Irr}( \SL_2(\cO_{2 \ell}) \mid \psi_{[A]})|=| \mathrm{Irr}(C_{\SL_2(\cO_{2\ell})}(\psi_{A}) \mid \psi_{[A]})|=|C_{\SL_2(\cO_{\ell})}(A)|=(q+1)q^{\ell-1}.$$

\end{proof}

\begin{proposition}\label{prop:second degrre for even r}
	 Let $\cO \in \dvrtwozero$ with ramification index $\ee $ and $\cO' \in \dvrtwoplus$ such that $\cO/\wp \cong \cO'/\wp'.$  Then for any $\ell \geq \mathrm{e}+1,$  the coefficients of $X^{(q-1)q^{2\ell-1}}$ in $	\calp_{\SL_2(\cO_{2 \ell})}(X)$ and in $	\calp_{\SL_2(\cO'_{2 \ell})}(X)$ are not equal. 
\end{proposition}

\begin{proof}
	Note that  $
	\calp_{\SL_2(\cO_{2 \ell})}(X)   = \calp^\pr_{\SL_2(\cO_{2 \ell})}(X) + \calp_{\SL_2(\cO_{2\ell-1})} (X).$ %Similarly for $\cO'.$ 
	 By Lemma~\ref{lem:second max dim from cyclic}, we have that 
	  the coefficient of $X^{(q-1)q^{2\ell-1}}$ in  $\calp_{\SL_2(\cO_{2 \ell-1})}(X)$  is zero.
	 Therefore the coefficients of $X^{(q-1)q^{2\ell-1}}$ in  $\calp_{\SL_2(\cO_{2 \ell})}(X)$  and $\calp_{\SL_2(\cO'_{2 \ell})}(X)$ are  equal to those  in  $\calp^\pr_{\SL_2(\cO_{2 \ell})}(X)$  and $\calp^\pr_{\SL_2(\cO'_{2 \ell})}(X)$ respectively. By Lemma~\ref{lem:second high dim iff irreducible},  the coefficient of $X^{(q-1)q^{2\ell-1}}$ in  $\calp^\pr_{\SL_2(\cO_{2 \ell})}(X)$ is 
	$\frac{(q^2-1) q^{2 \ell -2 }}{2}.$ We prove the proposition by showing that  the coefficient of $X^{(q-1)q^{2\ell-1}}$ in 
	$\calp^\pr_{\SL_2(\cO'_{2 \ell}) }(X)$ is strictly less than $\frac{(q^2-1) q^{2 \ell -2}}{2}.$ 
	For this, we focus our attention on representations $\rho \in \mathrm{Irr}(\SL_2(\cO'_{2 \ell})\mid \psi_{[A]})$ such that 
	%$\rho $ lies above $\psi_{[A]}$ and 
	$\bar{A}$ is irreducible. By Lemma~\ref{lem:dim divisibility}, these are the only representations that contribute, if at all, non-trivially to the coefficient of $X^{(q-1)q^{2\ell-1}}$ in 
	$\calp^\pr_{\SL_2(\cO'_{2 \ell}) }(X).$  Let $A = \mat 0 {a^{-1} \alpha} a \beta  \in M_2(\cO'_{\ell})$ such that $\bar{A}$ is irreducible, so in particular $\beta$ is invertible. By Lemma~\ref{lem:stabilizer-form}, we have  
	\[
	C_{\SL_2(\cO'_{2 \ell})} (\psi_{[A]}) = C_{\SL_2(\cO'_{2 \ell})} (\psi_{A})\{  I, \mat 1{\tilde{a}^{-1} \tilde \beta}01  \}. 
	\]

	 Note that $C_{\SL_2(\cO'_{2 \ell})} (\psi_{A})=C_S^{\ell}(\tilde{A})$ and $|C_{\SL_2(\cO'_{2 \ell})} (\psi_{A})|=|C_{\SL_2(\cO^\prime_{\ell})}(A)| \times |K^{\ell}|=\frac{|C_{\GL_2(\cO^\prime_{\ell})}(A)| \times |K^{\ell}|}{ |\det(C_{\GL_2(\cO'_{\ell})}(A))|}=(q+1)q^{4\ell-1}.$
	Here the last equality follows from  
	Lemmas~\ref{gl-centralizer-cardinality} and \ref{lem:image-det-map}.
	 So we have 
	 $\frac{|\SL_2(\cO'_{2\ell})|}{|C_S^{\ell}(\tilde{A})|}=(q-1)q^{2\ell-1}.$ By \cite[Lemma 3.2]{mypaper1}, every $\chi \in \mathrm{Irr}(C_S^\ell(\tilde{A})\mid \psi_{[A]})$ has dimension one.
	 Therefore, if  $\rho \in \mathrm{Irr}(\SL_2(\cO'_{2 \ell})\mid \psi_{[A]})$ 
	%$\rho$ is an irreducible representation of $\SL_2(\cO_{ 2\ell})$ lying above $\psi_{[A]}$ 
	such that $\dim(\rho ) = (q-1)q^{2\ell-1},$ then $\rho \cong \mathrm{Ind}_{C_S^\ell(\tilde{A})}^{\SL_2(\cO'_{2\ell})}\chi$ for some $\chi \in \mathrm{Irr}(C_S^\ell(\tilde{A})\mid \psi_{[A]}).$
	 Note that for such $\chi \in\mathrm{Irr}(C_S^\ell(\tilde{A})\mid \psi_{[A]}),$ 
	$\chi\neq \chi^{e_{\tilde{\beta}}}$ and 
	$ \mathrm{Ind}_{C_S^\ell(\tilde{A})}^{\SL_2(\cO'_{2\ell})}\chi \cong \mathrm{Ind}_{C_S^\ell(\tilde{A})}^{\SL_2(\cO'_{2\ell})}\chi^{e_{\tilde{\beta}}}.$ 
	%Hence the number of representations in  $\mathrm{Irr}(\SL_2(\cO'_{2 \ell}) \mid \psi_{[A]} )$ with dimension $(q+1)q^{2\ell-1}$ is equal to $\mathrm{n}_A :=  \frac{|\{ \chi \in \mathcal{E}_{\tilde{A}} \,\,\mid\,\, H_{\tilde{A}}(\chi)=C_{\SL_2(\cO'_{2\ell})}(\psi_A) \}| }{2}.$ 
	%
	%
	%
%
	Hence, for $\psi_{[A]} \in \Sigma_{\cO'}^{\IR},$ the number of  representations in  $\mathrm{Irr}(\SL_2(\cO'_{2 \ell}) \mid \psi_{[A]} )$ with dimension $(q-1)q^{2\ell-1}$ is equal to 
	 $ \mathrm{n}_A := \frac{|\{ \chi \in \mathrm{Irr}(C_S^\ell(\tilde{A})\mid \psi_{[A]}) \,\,\mid\,\,   H_{\tilde{A}}(\chi)=C_{\SL_2(\cO'_{2\ell})}(\psi_A) \}| }{2},$ where  $H_{\tilde{A}}(\chi)$ denote the stabilizer of $\chi$ in $C_{\SL_2(\cO'_{2\ell})}(\psi_{[A]}).$
	  Therefore the coefficient of $X^{(q-1)q^{2\ell-1}}$ in 
	$\calp^\pr_{\SL_2(\cO'_{2 \ell}) }(X)$  is equal to $\underset{ \psi_{[A]} \in \Sigma_{\cO'}^{\IR}}{\sum} \mathrm{n}_A,$ where   $\Sigma_{\cO'}^{\IR}$ is the set of orbits of $\psi_{[A]}$ with $\bar{A}$  irreducible.  
	For every $\psi_{[A]} \in \Sigma_{\cO'}^{\IR},$ we have the following 
	\[ 
	\mathrm{n}_A \leq \frac{|\mathrm{Irr}(C_S^\ell(\tilde{A})\mid \psi_{[A]})|}{2} = \frac{|C_{\SL_2(\cO_{\ell})}(A)|}{2}= \frac{(q+1)q^{\ell-1}}{2}. 
	\]
	
	Further by \cite[Proposition~7.5]{mypaper1} the inequality $\mathrm{n}_A \leq \frac{|\mathrm{Irr}(C_S^\ell(\tilde{A})\mid \psi_{[A]})|}{2}$ is strict for any $A \in \Sigma_{\cO'}^{\IR}$ such that $\mathrm{trace}(A) \in ( (\cO'_{\ell})^\times )^2.$ The existence of such a irreducible $A$ is easy to see, for example  $A=\mat{0}{1}{1}{1}$  will do. Combining this with $|\Sigma_{\cO'}^{\IR}|=(q-1)q^{\ell-1} $ (see \cite[Theorem~2.1]{mypaper1}),
	 %Lemma~\ref{lem:number-of-orbits-split semisimple}, 
	 we obtain that the  coefficient of $X^{(q-1)q^{2\ell-1}}$ in 
	$\calp^\pr_{\SL_2(\cO'_{2 \ell}) }(X)$  is strictly less than $\frac{(q^2-1) q^{2 \ell -2}}{2}.$ This completes the proof of the proposition. 
	
\end{proof}

\begin{proof}[Proof of Theorem~\ref{thm:zeta poly not eaqal- gp algebra}]
 For a finite group G, let $\Delta_m(G)$ be the coefficient of $X^m$ in  $\calp_{G }(X).$ For each $\ell' \geq \mathrm{e}+1,$ from Propositions~\ref{prop:second degrre for odd r} and \ref{prop:second degrre for even r} 
  we have  
 $  \Delta_{(q-1)q^{2\ell'-1}}(\SL_2(\cO_{2\ell'+1}))=\Delta_{(q-1)q^{2\ell'-1}}(\SL_2(\cO_{2\ell'})),$  
 $ \Delta_{(q-1)q^{2\ell'+1}}(\SL_2(\cO'_{2\ell'+1}))=\Delta_{(q-1)q^{2\ell'-1}}(\SL_2(\cO'_{2\ell'}))$
and 
 $\Delta_{(q-1)q^{2\ell'-1}}(\SL_2(\cO_{2\ell'}))\neq \Delta_{(q-1)q^{2\ell'-1}}(\SL_2(\cO'_{2\ell'})).$ 
  Therefore for $r\geq 2 \mathrm{e}+2,$ we obtain that  $\Delta_{(q-1)q^{2\lfloor\frac{r}{2} \rfloor-1}}(\SL_2(\cO_{r})) \neq \Delta_{(q-1)q^{2\lfloor\frac{r}{2} \rfloor-1}}(\SL_2(\cO'_{r})).$ This completes the proof.

\end{proof}

\begin{remark}\label{rmk:group algebra even r}
	Assume $\cO$ and $\cO'$ as in Theorem~\ref{thm:zeta poly not eaqal- gp algebra}. Then for any $\ell' > \mathrm{e},$ %Hassain and Singla 
	we showed in \cite{mypaper1} that both $\calp_{\SL_2(\cO_{2 \ell'})}(X)$ and $\calp_{\SL_2(\cO'_{2 \ell'}) }(X)$ are polynomials of degree $(q+1) q^{2 \ell'-1}$ and  their leading coefficients  are not equal (see \cite[Section~7.2]{mypaper1}).
	% This is an alternative proof of Theorem~\ref{thm:zeta poly not eaqal- gp algebra} for even $r.$
\end{remark}

\begin{remark}
	Assume $\cO$ and $\cO'$ as in Theorem~\ref{thm:zeta poly not eaqal- gp algebra}. Then for any $\ell' > \mathrm{e},$ 
	by combining Lemmas~\ref{lem:odd r dims from cyclic} and \ref{lem:max dim from cyclic}, we obtain that both $\calp_{\SL_2(\cO_{2 \ell'+1})}(X)$ and $\calp_{\SL_2(\cO'_{2 \ell'+1}) }(X)$ are polynomials of degree $(q+1) q^{2 \ell'}.$ Further by direct computation,
	%Theorem~\ref{thm:construction-SL-r odd- char=0}
	%
	it is easy to observe that the leading coefficients of  $\calp_{\SL_2(\cO_{2 \ell'+1})}(X)$ and $\calp_{\SL_2(\cO'_{2 \ell'+1}) }(X)$ are equal.
\end{remark}

Next we prove the following proposition. Which, along with Remark~\ref{rmk:group algebra even r},  give an alternate proof of Theorem~\ref{thm:zeta poly not eaqal- gp algebra} for the case $q\neq 2.$
\begin{proposition}
	Let $\cO \in \dvrtwozero$ with ramification index $\ee $ and $\cO' \in \dvrtwoplus$ such that $\cO/\wp \cong \cO'/\wp'.$ 
		%Let $\cO$ and $\cO'$ be compact discrete valuation rings  such that $\cO/\wp \cong \cO'/\wp',$ $2 \mid |\cO/\wp|,$ $\Char(\cO) =0$ with ramification index $\mathrm{e},$  and $\Char(\cO') = 2.$
	 Assume $|\cO/\wp|\neq2.$ Then for any $\ell' > \mathrm{e},$ the coefficients of $X^{(q+1)q^{2\ell'-1}}$ in $	\calp_{\SL_2(\cO_{2 \ell'+1})}(X)$ and in $	\calp_{\SL_2(\cO'_{2 \ell'+1})}(X)$ are not equal.
	 
\end{proposition}
\begin{proof}
	Recall that 
	$
	\calp_{\SL_2(\cO_{2 \ell'+1})}(X)   = \calp^\pr_{\SL_2(\cO_{2 \ell'+1})}(X) + \calp_{\SL_2(\cO_{2\ell'})} (X).
	$
	Since $\ell' > \mathrm{e},$ by Remark~\ref{rmk:group algebra even r} we have the coefficients of $X^{(q+1)q^{2\ell'-1}}$ in $\calp_{\SL_2(\cO_{2\ell'})} (X)$ and $\calp_{\SL_2(\cO'_{2\ell'})} (X)$ are not equal. Therefore to show  the proposition, 
	it is enough to show that   the coefficients of $X^{(q+1)q^{2\ell'-1}}$ in $\calp^\pr_{\SL_2(\cO_{2\ell'+1})} (X)$ and in $\calp^\pr_{\SL_2(\cO'_{2\ell'+1})} (X)$ are zero.
	For $\rho \in \mathrm{Irr}^\pr(\SL_2(\cO_{2 \ell'+1})),$ 
	by Lemma~\ref{lem:odd r dims from cyclic}, we have $\dim(\rho)$ is either in $\{(q+1)q^{2\ell'},(q-1)q^{2\ell'}\}$ or divisible by $(q^2 -1).$ Therefore if $\dim(\rho)=(q+1)q^{2\ell'-1},$ then $(q^2 -1)$ must divide $(q+1)q^{2\ell'-1}.$ This is not possible because $q^2-1\neq q+1$ for $q\neq 2.$ Hence the coefficients of $X^{(q+1)q^{2\ell'-1}}$ in $\calp^\pr_{\SL_2(\cO_{2\ell'+1})} (X)$ and $\calp^\pr_{\SL_2(\cO'_{2\ell'+1})} (X)$ are zero. Hence the  proposition.

\end{proof}

\bibliography{refs}{}
\bibliographystyle{siam}

\end{document}